\documentclass{amsart}

\usepackage{fullpage}
\usepackage{amsmath}
\usepackage{amsfonts}
\usepackage{amssymb}
\usepackage{yhmath}

\input xy
\xyoption{matrix}\xyoption{arrow}\xyoption{curve}  \xyoption{frame}

\begin{document}


\newcommand{\Hom}{\mathrm{Hom}}
\newcommand{\RHom}{\mathrm{RHom}^*}
\newcommand{\HOM}{\mathrm{HOM}}
\newcommand{\stHom}{\underline{\mathrm{Hom}}}
\newcommand{\Ext}{\mathrm{Ext}}
\newcommand{\Tor}{\mathrm{Tor}}
\newcommand{\HH}{\mathrm{HH}}
\newcommand{\Endo}{\mathrm{End}}
\newcommand{\ENDO}{\mathrm{END}}
\newcommand{\stEndo}{\mathrm{\underline{End}}}
\newcommand{\Tr}{\mathrm{Tr}}
\newcommand{\res}{\mathrm{res}}


\newcommand{\coker}{\mathrm{coker}}
\newcommand{\aut}{\mathrm{Aut}}
\newcommand{\op}{\mathrm{op}}
\newcommand{\ob}{\mathrm{ob}}

\newcommand{\add}{\mathrm{add}}
\newcommand{\ADD}{\mathrm{ADD}}
\newcommand{\ind}{\mathrm{ind}}
\newcommand{\rad}{\mathrm{rad}}
\newcommand{\soc}{\mathrm{soc}}
\newcommand{\ann}{\mathrm{ann}}
\newcommand{\im}{\mathrm{im}}
\newcommand{\chr}{\mathrm{char}}
\newcommand{\pdim}{\mathrm{p.dim}}
\newcommand{\cone}{\mathrm{Cone}}
\newcommand{\rk}{\mathrm{rk}}


\newcommand{\rmod}{\mbox{mod-}}
\newcommand{\Rmod}{\mbox{Mod-}}
\newcommand{\lmod}{\mbox{-mod}}
\newcommand{\lMod}{\mbox{-Mod}}
\newcommand{\stmod}{\mbox{\underline{mod}-}}
\newcommand{\stlmod}{\mbox{-\underline{mod}}}

\newcommand{\gmod}[1]{\mbox{mod}_{#1}\mbox{-}}
\newcommand{\gMod}[1]{\mbox{Mod}_{#1}\mbox{-}}
\newcommand{\Bimod}[1]{\mathrm{Bimod}_{#1}\mbox{-}}

\newcommand{\lrp}{\mbox{lrp}}
\newcommand{\stlrp}{\mbox{\underline{lrp}}}

\newcommand{\proj}{\mbox{proj-}}
\newcommand{\lproj}{\mbox{-proj}}
\newcommand{\Proj}{\mbox{Proj-}}
\newcommand{\inj}{\mbox{inj-}}
\newcommand{\coh}{\mbox{coh-}}


\newcommand{\und}[1]{\underline{#1}}
\newcommand{\gen}[1]{\langle #1 \rangle}
\newcommand{\floor}[1]{\lfloor #1 \rfloor}
\newcommand{\ceil}[1]{\lceil #1 \rceil}
\newcommand{\bnc}[2]{ \left( \begin{smallmatrix} #1 \\ #2 \end{smallmatrix}  \right)}

\newcommand{\bimo}[1]{{}_{#1}#1_{#1}}
\newcommand{\ses}[5]{\ensuremath{0 \rightarrow #1 \stackrel{#4}{\longrightarrow} 
#2 \stackrel{#5}{\longrightarrow} #3 \rightarrow 0}}
\newcommand{\bfa}{\mathbf{a}}
\newcommand{\bfb}{\mathbf{b}}
\newcommand{\bfc}{\mathbf{c}}
\newcommand{\bfd}{\mathbf{d}}
\newcommand{\A}{\mathcal{A}}
\newcommand{\B}{\mathcal{B}}
\newcommand{\C}{\mathcal{C}}
\newcommand{\D}{\mathcal{D}}
\newcommand{\E}{\mathcal{E}}
\newcommand{\F}{\mathcal{F}}
\newcommand{\N}{\mathcal{N}}
\newcommand{\T}{\mathcal{T}}
\newcommand{\s}{\mathcal{S}}
\newcommand{\X}{\mathcal{X}}
\newcommand{\Y}{\mathcal{Y}}
\newcommand{\K}{\mathcal{K}}
\newcommand{\ppX}{{}^{\perp}\mathcal{X}^{\perp}}
\newcommand{\tC}{\tilde{\mathcal{C}}}
\newcommand{\tK}{\tilde{K}(A)}
\newcommand{\ul}[1]{\underline{#1}}


\newtheorem{therm}{Theorem}[section]
\newtheorem{defin}[therm]{Definition}
\newtheorem{propos}[therm]{Proposition}
\newtheorem{lemma}[therm]{Lemma}
\newtheorem{coro}[therm]{Corollary}
\newtheorem{rmk}[therm]{Remark}

\title{Stable auto-equivalences for local symmetric algebras}
\author{Alex Dugas}
\address{Department of Mathematics, University of the Pacific, 3601 Pacific Ave, Stockton CA 95211, USA}
\email{adugas@pacific.edu}


\begin{abstract} 
We construct nontrivial auto-equivalences of stable module categories for elementary, local symmetric algebras over a field $k$.  These auto-equivalences are modeled after the spherical twists of Seidel and Thomas and the $\mathbb{P}^n$-twists of Huybrechts and Thomas, which yield auto-equivalences of the derived category of coherent sheaves on a variety.  For group algebras of $p$-groups in characteristic $p$ we recover many of the auto-equivalences corresponding to endo-trivial modules.  We also obtain analogous auto-equivalences for local algebras of dihedral and semi-dihedral type, which are not group algebras.  

\end{abstract}

\maketitle

 \section{Introduction} 
\setcounter{equation}{0}

One goal of the representation theory of self-injective algebras is to classify such algebras (over an algebraically closed field $k$) up to stable or derived equivalence.  Rickard's Morita theory for derived categories \cite{MTDC} 
 has been applied to describe the derived equivalence classes of self-injective algebras in numberous cases, including algebras of finite representation type, as well as various classes of tame symmetric algebras (see, for instance, \cite{SAFTT}, for an overview).  
 On the other hand, identifying the stable equivalence classes of self-injective algebras has been a harder problem, as there is no comparable Morita theory for describing or constructing equivalences of stable module categories.  Many of the known examples of stable equivalences between self-injective algebras actually come from derived equivalences.  Indeed, Rickard has shown that any standard derived equivalence between self-injective algebras induces a stable equivalence of Morita type \cite{DCSE}.  In the other direction, motivated in large part by Brou\'e's Abelian defect group conjecture, many authors have investigated the problem of `lifting' a stable equivalence to a derived equivalence (see, eg., \cite{Oku, EDCSA, Asa2}).   Although it is well-known that not every stable equivalence can be lifted to a derived equivalence in this way, relatively few examples of such non-liftable stable equivalences are available.  For self-injective algebras of finite representation type, for instance, all stable equivalences are induced by standard derived equivalences (up to natural isomorphism) \cite{Asa2, MutSym, CKL}.   In fact the majority of the known examples of stable equivalences that do not lift to derived equivalences arise in modular representation theory: for example, between blocks of group algebras with nonabelian Sylow $p$-subgroups having the trivial intersection property, as well as for group algebras of certain $p$-groups.

In order to better understand these examples and realize new ones, we focus here on constructing non-trivial auto-equivalences of $\stmod A$ for local symmetric $k$-algebras $A$.  Local algebras provide a natural starting point as they do not admit any nontrivial derived equivalences.   In particular, once we have a stable equivalence between local algebras it is relatively straightforward to check if it lifts to an equivalence of derived categories.

 The inspiration behind our constructions comes from the theory of endo-trivial modules for $p$-groups.  Recall that for a finite group $G$, and a field $k$ of characteristic $p$ dividing the order of $G$, a $kG$-module $M$ is said to be {\it endo-trivial} if $\Endo_{k}(M) \cong k \oplus P$ as $kG$-modules, where $k$ is the trivial $kG$-module and $P$ is projective.  (In fact, for a $p$-group $G$, this is equivalent to $\stEndo_{kG}(M) \cong k$ by a result of Carlson \cite{Car1}.) In this case the adjoint pair of exact functors $-\otimes_k M$ and $-\otimes_k M^*$ induces quasi-inverse auto-equivalences of $\stmod kG$.  Among these auto-equivalences are all (co)syzygy functors $\Omega^n$ for $n \in \mathbb{Z}$, corresponding to the endo-trivial modules $\Omega^n(k)$.  More interesting examples of endo-trivial modules have been exhibited for generalized dihedral, semidihedral and quaternion groups, and more generally for  $p$-groups with a maximal elementary abelian subgroup of rank $2$ (see \cite{CarThe2}).  
For these endo-trivial modules, the associated auto-equivalences of $\stmod kG$ are not induced by any auto-equivalences of the derived category $D^b(\rmod kG)$.

However, the existence of these auto-equivalences of $\stmod kG$ appears to depend on the Hopf algebra structure of $kG$ to ensure that the tensor product of two $kG$-modules is a $kG$-module.   We are interested in describing these auto-equivalences in terms only of the triangulated structure of the stable category in the hope that they can be extended to more general symmetric algebras.  In this direction, Alperin's construction of endo-trivial modules as relative syzygies \cite{Alp} affords some clues as to how these auto-equivalences can be generalized via a cone construction in the stable category.

  To summarize our main results, we let $A$ be a local, symmetric $k$-algebra, which is free as both a left and a right module over a subalgebra $R \cong k[t]/(t^m)$.  We define an $(A,A)$-bimodule $K$ as the kernel of the multiplication map $\mu : A \otimes_R A \rightarrow A$.  Our first construction provides an easily checked sufficient condition for $-\otimes_A K$ to induce an auto-equivalence of $\stmod A$.  

\begin{therm}[cf. Theorem 5.4]  Let $A, R$ and $K$ be as above, and assume that $\stEndo_A(k \otimes_R A) \cong k[\psi]/(\psi^2)$.  Then $-\otimes_A K$ induces an exact auto-equivalence of $\stmod A$.
\end{therm}

In this case, we call the auto-equivalence $\sigma_R := -\otimes_A \Omega^{-1}_{A^e}(K)$ a {\it spherical stable twist}.  Notice that $\Omega^{-1}_{A^e}(K) \cong \cone(\mu)$ in $\stmod A^e$, and thus for a typical $X$ in $\stmod A$, $\sigma_R (X) \cong \cone(X \otimes_R A \stackrel{\mu_X}{\longrightarrow} X)$ in $\stmod A$   Moreover, we show that $\sigma_R$ does not agree with any power of the syzygy functor, and hence is not induced by any auto-equivalence of $D^b(\rmod A)$.  We give some examples in Section 7 showing that these auto-equivalences are indeed more general than those given by endo-trivial modules for group algebras.  In fact, to our knowledge, this construction provides the first examples of nontrivial auto-equivalences of $\stmod A$ that do not come from derived equivalences when $A$ is not a group algebra.

There is an interesting parallel between these spherical stable twists and the spherical twists constructed on derived categories by Seidel and Thomas \cite{SeTh}. That is not to say that these auto-equivalences {\it are} spherical twists: in general one does not expect stable module categories of self-injective algebras to contain spherical objects as defined in \cite{SeTh} (see the discussion in Section 5).  However the condition $\stEndo_A(k \otimes_R A) \cong k[\psi]/(\psi^2)$, which is shared by $0$-spherical objects, is crucial in our proof.  Moreover, if we relax this condition to $\stEndo_A(k \otimes_R A) \cong k[\psi]/(\psi^{n+1})$ for $n \geq 1$, we can again construct auto-equivalences of $\stmod A$, but now using a double cone construction -- completely analogous to the $\mathbb{P}^n$ twists of Huybrechts and Thomas \cite{HuTh}.  

\begin{therm}[cf. Theorem 6.2] Let $A$ and $R$ be as above, and assume that $\stEndo_A(k \otimes_R A) \cong k[\psi]/(\psi^{n+1})$ for some $n \geq 1$, with $\psi$ induced by multiplication by some $y \in A$ that centralizes $R$.  If we set $$Q \cong \cone ( \cone (A \otimes_R A \stackrel{y\otimes 1 - 1 \otimes y}{\longrightarrow} A \otimes_R A) \stackrel{\bar{\mu}}{\longrightarrow} A)$$ in $\stmod A^e$, then $-\otimes_A Q$ induces an exact auto-equivalence of $\stmod A$.
\end{therm}

In this case we call the auto-equivalence $\rho_R = -\otimes_A Q$ a {\it $\mathbb{P}^n$ stable twist}.  Unlike the spherical stable twists, these auto-equivalences are not automatically nontrivial.  We show that if $n=1$ they correspond to the square of a spherical stable twist.  For larger $n$, the only examples we are presently aware of are actually over group algebras of $p$-groups.  Nevertheless in Section 7 we show that the auto-equivalences of $\stmod kG$ which are given by endo-trivial modules for certain $p$-groups $G$ are in fact of this form.  

Our proofs apply a strategy of Bridgeland's to show that an exact functor is an equivalence of triangulated categories.  While these ideas have been in use for some time in algebraic geometry to verify that a given Fourier Mukai transform induces an equivalence of derived categories (see \cite{FMTAG}), they have received little attention in the representation theory of algebras.  One exception is J. Grant's recent work constructing auto-equivalences of the derived category of a symmetric algebra, which he terms periodic twists \cite{Gra}.   As our focus is on stable module categories, we prove a version of Bridgeland's results with respect to a stronger notion of a spanning class, which is better suited for our context (see Proposition 2.2(2)), and which leads to the following nice criterion for an exact functor $F : \rmod A \rightarrow \rmod B$ to induce a stable equivalence.

\begin{therm}[cf. Theorem 3.4] Suppose $A$ and $B$ are split, indecomposable symmetric $k$-algebras with $B$ of Loewy length greater than $2$, and let $S_1, \ldots, S_n$ be the simple $A$-modules (up to isomorphism).  An exact functor $F : \stmod A \rightarrow \stmod B$ is an equivalence if and only if \begin{enumerate}
\item $\displaystyle \stHom_B(F(S_i), F(S_j)) \cong \left\{ \begin{array}{rl} k, & \mbox{if}\ i=j \\ 0, & \mbox{if}\ i\neq j \end{array} \right.$; and
\item $F$ induces monomorphisms $\displaystyle \Ext^1_A(S_i, S_j) \rightarrow \Ext^1_B(F(S_i),F(S_j))$ for all $i,j$.
\end{enumerate}
 \end{therm}

The paper is organized as follows.  We start off in Section 2 by reviewing, and modifying, some results of Bridgeland's that give sufficient conditions for an exact functor between triangulated categories to be fully faithful or an equivalence.  In Section 3 we specialize to stable module categories for symmetric algebras, and the exact functors between such categories that are induced by tensoring with so-called {\it left-right projective} bimodules.  Then in section 4 we specialize further to the situation where we have an elementary, local symmetric $k$-algebra $A$ that is free as a module on both sides over a subalgebra $R$ isomorphic to $k[t]/(t^m)$.  We prove some preliminary results under the assumption that the induced module $k\otimes_R A$ has a stable endomorphism algebra isomorphic to $k[\psi]/(\psi^{n+1})$.  The case where $n=1$ is considered in greater detail in Section 5, where we introduce the spherical stable twists mentioned above, and prove that they give non-trivial auto-equivalences of $\stmod A$.  In Section 6, we define  $\mathbb{P}^n$ stable twists (for $n \geq 1$) and show that they also yield auto-equivalences of $\stmod A$.  Finally, Section 7 returns to some of the examples of endo-trivial modules for $p$-groups.  We show that these endo-trivial modules are obtained as the images of the trivial module under a spherical or $\mathbb{P}^n$ stable twist, and thus these twists provide a description of the corresponding auto-equivalences of $\stmod kG$ that do not rely on the Hopf algebra structure of $kG$.  Moreover, it follows that these auto-equivalences continue to be defined for local algebras of dihedral or semidihedral type, which are not group algebras.  We conclude with some open questions.
 
Throughout this article, $k$ denotes a field and all categories and functors that we consider are $k$-linear.  For a $k$-algebra $\Lambda$, we write $\rmod \Lambda$ for the category of finitely presented right $\Lambda$-modules and $\stmod \Lambda$ for the associated stable category obtained by factoring out morphisms that factor through a projective.  As we typically work with right modules, we shalll write morphisms on the left and compose from right to left.  We also follow this convention for morphisms in abstract categories.  For a category $\C$, we shall write $\C(X,Y)$ for the set of morphisms from $X$ to $Y$ in $\C$.

\section{Equivalences of triangulated categories}
\setcounter{equation}{0}

We begin by reviewing some general results that are useful in establishing that a given exact functor between triangulated categories is an equivalence.  Let $k$ be a field and let $\T$ be a Hom-finite triangulated $k$-category with suspension $\Sigma$.  For objects $X, Y \in \T$ we write $\T(X,Y)$ for the morphisms from $X$ to $Y$, and we occasionally write $X[n]$ for $\Sigma^nX$.  Recall that an exact auto-equivalence $\nu$ of $\T$ is said to be a Serre functor if there are natural isomorphisms $\T(X,Y) \cong D\T(Y,\nu X)$ for all $X, Y \in \T$. 

\begin{defin} A collection $\C$ of objects in a triangulated category $\T$ is a {\bf spanning class} of $\T$ if \begin{itemize} \item $\T(X,Y[i]) = 0$ for all $X \in \C$ and all $i \in \mathbb{Z}$ implies $Y \cong 0$; and \item $\T(Y[i],X) = 0$ for all $X \in \C$ and all $i\in \mathbb{Z}$ implies $Y \cong 0$. \end{itemize}  A collection $\D$ of objects in $\T$ will be said to be a {\bf strong spanning class} if \begin{itemize} \item $\T(X,Y) = 0$ for all $X \in \D$ implies $Y \cong 0$; and \item $\T(Y,X) = 0$ for all $X \in \D$ implies $Y \cong 0$.\end{itemize}
\end{defin}

Observe that if $\T$ has a Serre functor, then the two conditions in either of the above definitions are equivalent.  Thus it suffices to check only one of the conditions in this case.  In particular, when $\T$ has a Serre functor, useful examples of (strong) spanning classes are of the form $\C = \{X\} \cup X^{\perp}$ for any $X \in \T$, where $X^{\perp} = \{Y \in \T\ |\ \T(X,Y)=0\}$.  Another example of a strong spanning class is a maximal system of orthogonal bricks (see \cite{SMS}).

Now suppose that $\T$ and $\T'$ are triangulated categories with suspensions $\Sigma$ and $\Sigma'$ respectively.  If $F: \T \rightarrow \T'$ is an exact equivalence of triangulated categories, then clearly it has a left and right adjoint, and it can be seen that this adjoint is also exact.  More generally, suppose that $(F,\sigma) : \T \rightarrow \T'$ is an exact functor, where $\sigma : F\Sigma \stackrel{\cong}{\longrightarrow} \Sigma' F$ is an isomorphism, and suppose that $H : \T' \rightarrow \T$ is right adjoint to $F$.  Notice that this means that $H$ is a $k$-linear functor and there is a natural isomorphism $\eta : \T'(F(-),-) \stackrel{\cong}{\longrightarrow} \T(-,H(-))$.  We can define an isomorphism $\tau : H\Sigma' \rightarrow \Sigma H$ via the sequence of isomorphisms of bi-functors 
\begin{equation}
 \vcenter{\xymatrixcolsep{4.0pc} \xymatrix{\T(\Sigma(-),H\Sigma'(-)) \ar[r]^{\eta^{-1}} \ar@{-->}[d]^{\T(\Sigma(-),\tau)} &  \T'(F\Sigma(-),\Sigma'(-)) \ar[r]^{\T'(\sigma, \Sigma'(-))^{-1}} &  \T'(\Sigma'F(-), \Sigma'(-)) \ar[d]^{\Sigma'^{-1}} \\ \T(\Sigma(-),\Sigma H(-)) & \T(-,H(-)) \ar[l]_{\Sigma} & \T'(F(-),-). \ar[l]_{\eta} }}
\end{equation}
We shall write $h : 1 \rightarrow HF$ for the unit of the adjunction, so that $h_Y = \eta(1_{FY})$ for any object $Y$ of $\T$.  Starting with $1_{FY}$ in the bottom right of the above diagram, we see that

 \begin{equation} \Sigma(h_Y) = \Sigma(\eta(1_{FY})) = \tau_{FY} \circ \eta(\Sigma'(1_{FY})\sigma_Y) = \tau_{FY} \circ \eta(\sigma_Y) =  \tau_{FY} \circ H(\sigma_Y)\circ h_{\Sigma Y}. \end{equation}
 This equality expresses the fact that the unit of the adjunction is compatible with the suspension functor in the appropriate sense.  Dually, for the co-unit $g: FH \rightarrow 1$, one obtains $\Sigma'(g_X) \circ \sigma_{HX} \circ F(\tau_X) = g_{\Sigma' X}$ for any object $X$ of $\T'$. 

\begin{propos}[Cf. \cite{Bri}] Let $F : \T \rightarrow \T'$ be an exact functor between triangulated categories with a right adjoint $H$ and a left adjoint $G$.  Then $F$ is fully faithful if either of the following two conditions is satisfied.
\begin{enumerate}
\item There exists a spanning class $\C$ of $\T$ such that the natural homomorphisms $$F : \T(X,Y[i]) \rightarrow \T'(FX,F(Y[i]))$$ are bijective for all $X, Y \in \C$ and all $i \in \mathbb{Z}$. 
\item There exists a strong spanning class $\D$ of $\T$ such that for all $X, Y \in \D$ the natural homomorphisms $$F : \T(X,Y[i]) \rightarrow \T'(FX,F(Y[i]))$$ are bijective for $i = 0$ and injective for $i=1$.
\end{enumerate}
\end{propos}

\noindent
{\it Proof.}   The first condition corresponds to a theorem of Bridgeland in \cite{Bri, FMTAG}.  We prove the sufficiency of the second condition in a similar manner.  Consider the commutative diagram

\begin{equation}
\vcenter{\xymatrix{\T(X,Y) \ar[r]^{h_Y \circ} \ar[d]_{\circ g_X} \ar[dr]^F & \T(X,HFY) \ar[d]^{\cong} \\ \T(GFX,Y) \ar[r]^{\cong} & \T'(FX,FY)}}
\end{equation} 
where $g$ denotes the co-unit of the adjunction $G \dashv F$ and $h$ denotes the unit of the adjunction $F \dashv H$.  Let $X \in \D$ and construct a triangle $GFX \stackrel{g_X}{\longrightarrow} X \longrightarrow Z \rightarrow$.  If we apply $\T(-,Y)$ to this triangle we get the exact sequence in the top row of the commutative diagram.
$$\xymatrixcolsep{3.0pc} \xymatrix{\T(X,Y) \ar[r]^{\circ g_X} \ar[dr]_F & \T(GFX,Y) \ar[r] \ar[d]^{\cong} & \T(Z[-1],Y) \ar[r] & \T(X[-1],Y) \ar[r]^{\circ \Sigma^{-1}(g_X)}  \ar[d]^{\cong}_{\Sigma} & \T(GFX[-1],Y) \ar[d]^{\cong}_{\Sigma} \ar[r] & 
\\ & \T'(FX,FY) & & \T(X,Y[1]) \ar[r]^{\circ g_X} \ar[dr]_F & \T(GFX, Y[1]) \ar[d]^{\cong} \\
& & & & \T'(FX,F(Y[1]))}$$
If $Y \in \D$, the first map labeled $F$ is bijective and the second is injective.  Thus the commutativity of the diagram yields $\T(Z[-1],Y) = 0$.  Consequently, $Z \cong 0$ and $g_X$ is an isomorphism.  It follows from the commutativity of (2.3) that for any $X \in \D$, $h_Y \circ : \T(X,Y) \rightarrow \T(X,HFY)$ is bijective for all $Y \in \T$.  We now apply $\T(X,-)$ with $X \in \D$ to the triangle $Y \stackrel{h_Y}{\longrightarrow} HFY \longrightarrow W \rightarrow$ to get $\T(X,W) = 0$.  Here we are using the fact that $\T(X,\Sigma(h_Y))$ is also bijective, which follows from (2.2) since we know $\T(X,h_{\Sigma Y})$ is bijective.  Hence $W \cong 0$ and $h_Y$ is an isomorphism for all $Y \in \T$.  Returning to the commutative square (2.3) we see that $F$ must induce an isomorphism for all $X, Y \in \T$.  $\Box$ \\

Once we have established that a given functor $F$ is fully faithful, we can use the following result to verify that $F$ is an equivalence of triangulated categories.

\begin{propos}[Bridgeland \cite{Bri}]  Let $F : \T \rightarrow \T'$ be a fully faithful exact functor between triangulated categories such that $\T$ contains nonzero objects and $\T'$ is indecomposable.  Then $F$ is an equivalence of categories if and only if $F$ has a left adjoint $G$ and a right adjoint $H$ such that $H(Y) \cong 0$ implies $G(Y) \cong 0$ for any $Y \in \T'$.
\end{propos}

Notice, in particular, that if $G$ is both a left and right adjoint to $F$ in the above situation, then $F$ is an equivalence if and only if it is fully faithful.

\section{Left-right projective bimodules}
\setcounter{equation}{0}
We now specialize the ideas of the previous section to the setting of stable categories of self-injective $k$-algebras over a field $k$.  From now on we only consider finite-dimensional self-injective $k$-algebras, and finite-dimensional modules.  If $\Lambda$ is self-injective, we write $\rmod \Lambda$ for the category of finte-dimensional right $\Lambda$-modules, $\stmod \Lambda$ for the associated stable category, and $D^b(\rmod \Lambda)$ for the bounded derived category.  The stable category $\stmod \Lambda$ is always triangulated with the cosyzygy functor $\Omega^{-1}$ as the suspension.  Since the most natural examples of exact functors between stable module categories are those induced by tensoring with bimodules, we begin by reviewing some important features of the relevant bimodules.

Let $A$ and $B$ be basic, indecomposable and non-semisimple self-injective $k$-algebras.  We consider $(A,B)$-bimodules on which $k$ acts centrally, which we identify with right modules over $A^{\op} \otimes_k B$.  If $A=B$, this ring is also called the enveloping algebra of $A$ and is denoted $A^e$.  We let $D = \Hom_k(-,k)$ be the duality with respect to the ground field.  For an $(A,B)$-bimodule $M$, we get a $(B,A)$-bimodule $DM$.   We may also take the dual of $M$ with respect to either $A$ or $B$ on one side.  We set ${}^*M = \Hom_A(M,A)$, the dual of $M$ with respect to its left-module structure, and $M^* = \Hom_B(M,B)$ the dual with respect to the right-module structure.

We are interested in functors of the form $-\otimes_A M_B : \rmod A \rightarrow \rmod B$ for a bimodule ${}_A M_B$.  This functor is exact if and only if ${}_A M$ is projective, and it takes projective $A$-modules to projective $B$-modules if and only if $M_B$ is projective.  Following \cite{ZP1}, we call the bimodule ${}_A M_B$ {\bf left-right projective} if ${}_A M$ and $M_B$ are both projective.  For such bimodules the functor $-\otimes_A M_B$ induces an exact functor between triangulated categories $\stmod A \rightarrow \stmod B$.

  Observe that $\Lambda := A^{\op} \otimes_k B$ is self-injective, and thus $\rmod \Lambda$ is a Frobenius category.  The subcategory $\lrp(A,B)$ of left-right projective bimodules is an exact subcategory containing the projective bimodules, and hence it is also a Frobenius category.  We write $\stlrp(A,B)$ for the corresponding stable category of left-right projective bimodules, and note that it inherits the structure of a triangulated category.  A sequence $$X \stackrel{\ul{f}}{\longrightarrow} Y \stackrel{\ul{g}}{\longrightarrow} Z \stackrel{\ul{h}}{\longrightarrow} X[1]$$ is a distinguished triangle if and only if there is a short exact sequence $\ses{X}{Y\oplus P}{Z}{\bnc{f}{p}}{(g\ q)}$ corresponding to the map $\ul{h} \in \Ext^1_\Lambda(Z,X) \cong \stHom_{\Lambda}(Z,\Omega^{-1}X)$.
  
\begin{lemma} Let $M$ be a left-right projective $(A,B)$-bimodule and let $\Lambda = A^{\op} \otimes_k B$.  Then we have isomorphisms in $\stlrp(A,B)$ $$\Omega_{\Lambda}(M) \cong M \otimes_B \Omega_{B^e}(B) \cong \Omega_{A^e}(A) \otimes_A M.$$
\end{lemma}

\noindent
{\it Proof.}  We have an exact sequence $\ses{\Omega_{A^e}(A)}{P}{A}{}{}$ in $\lrp(A,A)$ for a projective $A^{e}$-module $P$.  If we tensor on the right with $M$, we obtain an exact sequence $\ses{\Omega_{A^e}(A)\otimes_A M}{P \otimes_A M}{M}{}{}$ in $\lrp(A,B)$.  Since $M_B$ is projective, $P \otimes_A M$ is a projective $(A,B)$-bimodule, and it follows that $\Omega_{A^e}(A)\otimes_A M \cong \Omega_{\Lambda}(M)$ up to projective summands.  The other isomorphism is proved similarly. $\Box$\\

Let $M$ be a left-right projective $(A,B)$-bimodule.  Since $A$ and $B$ are self-injective $DM, M^*$ and ${}^*M$ are all left-right projective $(B,A)$-modules.  The functor $-\otimes_A M_B$ has a right adjoint $\Hom_B(M,-) \cong -\otimes_B M^*$ and a left adjoint $-\otimes_A {}^*M$.  Furthermore, it is not hard to see that these functors induce adjoint pairs of functors between $\stmod A$ and $\stmod B$.

\begin{lemma} For an $(A,B)$-bimodule $M$, we have isomorphisms $DA \otimes_A M \cong D(^*M)$ and $M \otimes_B DB \cong D(M^*)$ of $(A,B)$-bimodules.  Hence $M^* \cong {}^*M$ if and only if $DA \otimes_A M \cong M \otimes_B DB$ as bimodules.  In particular, if $A$ and $B$ are symmetric algebras, then $M^* \cong {}^*M \cong DM$.
\end{lemma}

\noindent
{\it Proof.}  Using Hom-tensor adjunction, we have natural isomorphisms 
\begin{eqnarray*} D(DA \otimes_A M) = \Hom_k(DA \otimes_A M, k) & \cong & \Hom_A( {}_A M, \Hom_k(DA, k)) \\ & \cong & \Hom_A (_A M, A) = {}^*M. \end{eqnarray*}
By naturality, these are bimodule isomorphisms.  The isomorphism $D(M^*) \cong M \otimes_B DB$ is proved similarly.  Finally, if $A$ and $B$ are symmetric, then $DA \cong A$ and $DB \cong B$ as bimodules.  Thus $M \cong D(M^*) \cong D(^*M)$, and applying the duality $D$ gives the desired result. $\Box$ \\

We have the following consequence of Proposition 2.3.  The assumptions on $B$ guarantee that $\stmod B$ is indecomposable.

\begin{coro}  Let ${}_A M_B$ be a left-right projective bimodule and assume $A$ and $B$ are symmetric algebras with $B$ indecomposable and of Loewy length greater than $2$, and let $S_1, \ldots, S_n$ be the simple $A$-modules (up to isomorphism).  Then $-\otimes_A M_B : \stmod A \rightarrow \stmod B$ is an equivalence if and only if it is fully faithful.
\end{coro}

Using the obvious fact that the set of simple $A$-modules is a strong spanning class in $\stmod A$, we can apply this Corollary and Proposition 2.2(2) to give a characterization of the left-right projective bimodules that induce stable equivalences.  This result might be seen as a rough analogue of Bridgeland's Theorem 1.1 in \cite{Bri} describing which vector bundles induce Fourier-Mukai transforms.  For simplicity, we assume that $A$ and $B$ are split over the field $k$.  In particular, every simple module has endomorphism ring isomorphic to $k$.

\begin{therm}\label{thm:equiv-criterion} Suppose $A$ and $B$ are split, indecomposable symmetric $k$-algebras with $B$ of Loewy length greater than $2$, and let $S_1, \ldots, S_n$ be the simple $A$-modules (up to isomorphism).  An exact functor $F : \stmod A \rightarrow \stmod B$ is an equivalence if and only if 
\begin{enumerate}
\item $\displaystyle \stHom_B(F(S_i), F(S_j)) \cong \left\{ \begin{array}{rl} k, & \mbox{if}\ i=j \\ 0, & \mbox{if}\ i\neq j \end{array} \right.$; and
\item $F$ induces monomorphisms $\displaystyle \Ext^1_A(S_i, S_j) \rightarrow \Ext^1_B(F(S_i),F(S_j))$ for all $i,j$.
\end{enumerate}
 \end{therm}

  \begin{rmk}\label{rmk:preserves-nonsplit} \textup{In case $F$ is induced by a functor $-\otimes_A M$ for a left-right projective bimodule $M$, the second condition of the above theorem may be phrased as: $\ses{S_j \otimes_A M}{U \otimes_A M}{S_i \otimes_A M}{}{}$ is non-split in $\rmod B$ for every non-split extension $\ses{S_j}{U}{S_i}{}{}$ in $\rmod A$.}
  \end{rmk}

Under some mild assumption on the ground field or the semisimple quotients of the algebras $A$ and $B$, any equivalence of stable categories induced by a left-right projective bimodule is of Morita type.   Recall that a pair of bimodules ${}_A M_B $ and ${}_B N_A$ is said to induce a {\it stable equivalence of Morita type} between $A$ and $B$ if we have isomorphisms of $(A,A)$ and $(B,B)$-bimodules, respectively:
$$M \otimes_B N \cong A \oplus P \ \ \mbox{and}\ \ \ N \otimes_A M \cong B \oplus Q$$ for projective bimodules $P$ and $Q$.  A (semisimple) $k$-algebra $R$ is said to be {\it separable} if $R \otimes_k K$ is semisimple for any field $K$ containing $k$.  This holds automatically if $k$ is perfect, or if $R$ splits over $k$.  In particular, the following theorem applies to the elementary, local symmetric $k$-algebras that we study in the following sections.

\begin{therm}[Rickard \cite{AMRT}]  Assume $A/\rad A$ and $B/\rad B$ are separable $k$-algebras, and ${}_A M_B$ is an indecomposable, left-right projective bimodule such that $-\otimes_A M_B$ induces an equivalence $\stmod A \rightarrow \stmod B$.  Then $M$ and $N := M^*$ give a stable equivalence of Morita type between $A$ and $B$.
\end{therm}

The following result of Linckelmann gives a criterion for a stable equivalence of Morita type to lift to a Morita equivalence.  In practice, it provides a useful way of comparing different stable equivalences of Morita type.

\begin{propos}[\cite{ML1}]\label{prop:LinckelmannTheorem}  Suppose that $A$ and $B$ are self-injective $k$-algebras with no projective simple modules, and ${}_A M_B$ is a left-right projective bimodule for which $-\otimes_A M_B$ induces a stable equivalence $\stmod A \rightarrow \stmod B$.  Then $-\otimes_A M_B :\rmod A \rightarrow \rmod B$ induces a Morita equivalence if and only if $S \otimes_A M_B$ is simple for each simple $A$-module $S$.
\end{propos}

On the other hand, it is often much harder to determine if a stable equivalence of Morita type lifts to an equivalence of derived categories.  Here, we say that a stable equivalence of Morita type $\alpha : \stmod \Lambda \rightarrow \stmod \Gamma$ {\it lifts to a derived equivalence} if there is an equivalence of triangulated categories $\tilde{\alpha} : D^b(\rmod \Lambda) \rightarrow D^b(\rmod \Gamma)$ and a natural isomorphism $\pi_\Gamma \tilde{\alpha} \cong \alpha \pi_{\Lambda}$ with $\pi_{\Lambda}$ denoting the localization functor $D^b(\rmod \Lambda)  \rightarrow D^b(\rmod \Lambda)/K^b(\proj \Lambda) \approx \stmod \Lambda$, and similarly for $\pi_{\Gamma}$.  Indeed, this problem has received much attention for its relevance to Brou\'e's conjecture.  For local algebras, however, it has a simple answer. (We continue to assume that $A$ and $B$ are basic and self-injective.)

\begin{propos}\label{prop:trivial-dpic} Suppose $\alpha : \stmod A \rightarrow \stmod B$ is a stable equivalence of Morita type with $A$ local.  Then the following are equivalent:
\begin{enumerate}
\item $\alpha$ can be lifted to a derived equivalence $\tilde{\alpha} : D^b(\rmod A) \rightarrow D^b(\rmod B)$;
\item There exists a ring isomorphism $\phi: A \rightarrow B$ such that  $\phi^* (\alpha(k)) \cong \Omega^n k$  in $\stmod A$ for some $n \in \mathbb{Z}$, where $k$ denotes the unique simple $A$-module and $\phi^* : \rmod B \rightarrow \rmod A$ is the equivalence induced by $\phi$.
\item There exists a Morita equivalence $F :\rmod A \rightarrow \rmod B$  such that $\alpha \cong F \circ \Omega^n$ for $n \in \mathbb{Z}$.  
\end{enumerate}
\end{propos}

\noindent
{\it Proof.}  The following arguments are standard, but we include them for convenience of the reader.  $(1) \Rightarrow (2):$  The derived equivalence $\tilde{\alpha}$ is associated to a tilting complex $T \in K^b(\proj A)$ such that $\tilde{\alpha}(T) \cong B$.  It is easy to see that the only basic tilting complexes over a local algebra $A$ are of the form $T = A[n]$ for some $n  \in \mathbb{Z}$.  Thus $\tilde{\alpha} \circ \Sigma^n$ takes $A$ to $B$ and we let $\phi : A \rightarrow B$ be the ring isomorphism induced by this equivalence on endomorphism rings.  We now apply Proposition 7.1 in \cite{MTDC} to the composite $\phi^* \tilde{\alpha}\Sigma^n$ to conclude that $\phi^* \tilde{\alpha} (X[n]) \cong X$ for all $X$ in $\D^b(\rmod A)$.  In particular, setting $X = k[-n]$ gives the desired conclusion upon projecting down to $\stmod A$.

$(2) \Rightarrow (3):$ Since $\phi^*\alpha$ is a stable equivalence of Morita type it commutes with $\Omega$, and we see that $G := \phi^*\alpha \Omega^{-n}$ must  be a Morita auto-equivalence of $A$ by Proposition~\ref{prop:LinckelmannTheorem}.  Letting $F$ be the Morita equivalence $(\phi^{-1})^* \circ G$ now yields $\alpha \cong F \circ \Omega^n$.

$(3) \Rightarrow (1):$ Any Morita equivalence $F$ induces an equivalence of derived categories, while any power $\Omega^n$ of the suspension functor of $\stmod A$ clearly lifts to the corresponding power of the suspension functor of $D^b(\rmod A)$.  Thus $\alpha$ can be lifted to an equivalence of derived categories as well.  $\Box$ \\

Since any Morita auto-equivalence of $A$ preserves the isomorphism class of $\Omega^n k$ for all $n$, we obtain the following simple criterion for non-liftability of an auto-equivalence of $\stmod A$.

\begin{coro}\label{coro:no-lift} Let $A$ be a local self-injective algebra with unique simple module $k$, and $\alpha : \stmod A \rightarrow \stmod A$ a stable equivalence of Morita type.  If $\alpha$ satisfies $\alpha(k) \ncong \Omega^n(k)$ for all $n \in \mathbb{Z}$, then $\alpha$ does not lift to an auto-equivalence of $D^b(\rmod A)$.
\end{coro}

\section{Set up and preliminary results}
\setcounter{equation}{0}

We now focus on an elementary, local, symmetric $k$-algebra $A$.   Recall that $A$ elementary means that $A/\rad A \cong \prod_{i=1}^s k$ as $k$-algebras for some $s \geq 1$, and we must have $s=1$ since $A$ is local.  Furthermore, it follows that $A$ is a split $k$-algebra with a unique simple (right) module $k \cong A/\rad A$, up to isomorphism.  Let $x \in \rad\ A$ be an element with $x^m = 0$ but $x^{m-1} \neq 0$ for some $m \geq 2$, and set $R = k\cdot 1 + k\cdot x + \cdots + k\cdot x^{m-1} \cong k[t]/(t^m)$, the unital subalgebra of $A$ generated by $x$.  We shall assume that $A$ is free as an $R$-module on both sides.  Thus we have exact induction and restriction functors $-\otimes_R A : \rmod R \rightarrow \rmod A$ and $-\otimes_A A_R : \rmod A \rightarrow \rmod R$.  Since $R$ and $A$ are symmetric algebras, and ${}_A A_R$ and ${}_R A_A$ are left-right projective bimodules, induction and restriction are left and right adjoints of each other.  Furthermore, they induce mutually adjoint exact functors between $\stmod R$ and $\stmod A$.   When convenient we will write $F = -\otimes_R A$ for induction and $G=-\otimes_A A_R$ for restriction.  We also let $\eta : \stHom_A(FX,Y) \stackrel{\cong}{\longrightarrow} \stHom_R(X,GY)$  and $\theta : \stHom_R(GX,Y) \stackrel{\cong}{\longrightarrow} \stHom_A(X,FY)$ denote the adjugants.

We also write $k$ for the unique simple (right) $R$-module $R/xR$, and we set $T_A := k \otimes_R A \cong A/xA$.  Notice that $\Endo_A(T) \cong \{a \in A\ |\ ax \in xA\}/xA$ and $\stEndo_A(T) \cong \{a \in A\ |\ ax \in xA\}/(xA + \ann_l(x))$.  Our key assumption in constructing auto-equivalences of $\stmod A$ is that $$\stEndo_A(T) \cong k[\psi]/(\psi^{n+1})$$ for some $n \geq 1$.  By the above description of $\stEndo_A(T)$, we may assume that $\psi$ is induced by left multiplication by an element $y \in \rad\ A$, which we henceforth fix.

\begin{lemma}\label{lemma:RestrictionPsi} Assume that $\stEndo_A(T) = k[\psi]/(\psi^{n+1})$ with $\psi$ induced by left multiplication by an element $y \in A$.  Then $T_R \cong k^{n+1} \oplus R^l$ for some integer $l \geq 0$ and the restriction of $\psi$ coincides with the following $(n+1)\times(n+1)$ matrix (with respect to an appropriate basis)  $$ \left( \begin{array}{ccccc} 0 & 0 &  0 & \cdots & 0 \\ 1 & 0 & 0 & \cdots & 0 \\ 0 & 1 & 0 & \cdots & 0 \\ \vdots & & \ddots & & \vdots \\
0 & 0 & \cdots & 1  & 0 \end{array} \right) \in \Endo_R(k^{n+1}) \cong \stEndo_R(T).$$
\end{lemma}

\noindent
{\it Proof.}  We first decompose $T_R = U_R \oplus R^l$ where $U_R$ is the sum of all nonprojective direct summands of $T_R$.  Since induction is left adjoint to restriction we have an isomorphism $\eta : \stHom_A(T,T) \stackrel{\cong}{\longrightarrow} \stHom_R(k,T_R) = \stHom_R(k,U) = \Hom_R(k,U)$.  Thus $\soc\ U$ has dimension $n+1$, and $U$ must decompose as a direct sum of $n+1$ indecomposable $R$-modules (which are all uniserial).  Furthermore, $\stHom_R(k,T_R) = \Hom_R(k,U)$ is spanned by the set $$\{\eta(1_T), \eta(\psi), \ldots, \eta(\psi^{n})\} = \{\eta(1_T), y\eta(1_T), \ldots, y^{n}\eta(1_T)\}.$$  Consequently $\soc\ U$ is spanned by $\{y^i\eta(1_T)(1)\}_{0\leq i \leq n}$.   Observe that $\eta(1_T)$ corresponds to the map $k \rightarrow A/xA$ sending $1 \in k$ to $1 + xA \in A/xA$, which is a split monomorphism.  Furthermore, if $f : A/xA \rightarrow k$ is a splitting for $\eta(1_T)$, then $f \circ \psi = 0$ as the map $\psi$ corresponds to multiplication by an element $y \in \rad A$ on which $f$ vanishes.

Since induction is also right adjoint to restriction, we have an isomorphism $\theta : \stHom_A(T,T) \stackrel{\cong}{\longrightarrow} \stHom_R(T_R,k) = \stHom_R(U,k) = \Hom_R(U,k)$.  Hence $\Hom_R(U,k)$ is spanned by $$\{\theta(1_T), \theta(\psi), \ldots, \theta(\psi^{n})\} = \{\theta(1_T), \theta(1_T)y, \ldots, \theta(1_T)y^{n}\}.$$  As the map $f : T \rightarrow k$ above has $f y = 0$, we must have $f = \theta(1_T)y^{n}$ up to a scalar multiple.  As $f(\eta(1_T)(1)) = 1$, we have $\theta(1_T)y^{n-i} (y^i\eta(1_T)(1))=1$ for each $i$.  Consequently, the elements of the socle of $U$ are not in $\rad\ U$ and it follows that they generate $U$.  Therefore $U = \soc\ U \cong k^{n+1}$, and $y$ has the desired matrix with respect to the basis $\{y^i\eta(1_T)(1)\}_{0\leq i \leq n}$.  $\Box$\\

\begin{rmk}\label{rmk:2dim} The converse of the above lemma also appears to be true.  However, when $n=1$, it is enough to verify that $T_R  \cong k^2 \oplus R^l$.  For then $\stEndo_A(T)$ is a two-dimensional local ring and hence must be isomorphic to $k[\psi]/(\psi^2)$. \\
\end{rmk}

When $m=2$, the short exact sequence $\ses{k}{R}{k}{}{}$ yields an exact sequence of $A$-modules $\ses{T}{A}{T}{}{}$ which shows that $\Omega T \cong T$.  For $m > 2$, the projective presentation $\ses{k \longrightarrow R}{R}{k}{x\cdot}{}$ yields a projective presentation for $T$ of the form $\ses{T \longrightarrow A}{A}{T}{x\cdot}{}$, showing that $\Omega^2T \cong T$.  In this case we will need some additional information about the maps between $\Omega T$ and $T$.

\begin{lemma}\label{lemma:RestrictionXi}  Let $R, A$ and $T = k \otimes_R A$ be as above.  Then $\stHom_A(\Omega T,T)$ has a $k$-basis $\{\xi, \psi\xi, \ldots, \psi^{n}\xi\}$ where $\xi = f \otimes_R A$ for a nonzero map $f: \Omega k_R \rightarrow k_R$.   Furthermore, provided that either $n=1$ or that $x^{m-2}$ commutes with $y$, the restriction of $\xi$ is given by the diagonal matrix $$\begin{pmatrix}f &  & 0 \\ & \ddots & \\ 0 & & f \end{pmatrix} \in \Hom_R(\Omega k^{n+1},k^{n+1}) \cong \stHom_R(\Omega T_R, T_R)$$ with respect to appropriate bases of the projective-free summands of $\Omega T_R$ and $T_R$,
\end{lemma}

\noindent
{\it Proof.}  We have an isomorphism $\eta: \stHom_A(\Omega T,T) \stackrel{\cong}{\rightarrow} \stHom_R(\Omega k_R, T_R)$.  Using the same basis for the projective-free part of $T_R \cong k^{n+1}$ as in the last lemma, we have $$\eta(\xi) = \eta(F(f)) = \eta(1_T)\circ f = \begin{pmatrix} 1 \\ 0 \\ \vdots \\ 0 \end{pmatrix}f = \begin{pmatrix} f \\ 0 \\ \vdots \\ 0 \end{pmatrix},$$ and 
$$ \eta(\psi^i \xi) = G(\psi^i)\eta(\xi) = \begin{pmatrix}  0 & 0 &  0 & \cdots & 0 \\ 1 & 0 & 0 & \cdots & 0 \\ 0 & 1 & 0 & \cdots & 0 \\ \vdots & & \ddots & & \vdots \\
0 & 0 & \cdots & 1  & 0 \end{pmatrix}^i \begin{pmatrix} f \\ 0 \\ \vdots \\ 0 \end{pmatrix} = \begin{pmatrix} 0 \\ \vdots \\ f  \\ \vdots \\ 0 \end{pmatrix}$$ with $f$ in the $(i+1)^{th}$ row.   Thus $\eta$ sends $\{ \psi^i \xi\}_{0 \leq i \leq n}$ to a basis for $\stHom_R(\Omega k_R, T_R)$ and the first conclusion follows.  

For the second claim, observe that $(\Omega T)_R \cong \Omega (T_R) \cong \Omega( k_R^{n+1}) \cong (\Omega k_R)^{n+1}$ in $\stmod R$.  As in the proof of Lemma 4.1, we fix such a decomposition of $\Omega T_R$ so that $\eta(\Omega(\psi^i))$ is the identity map from $\Omega k_R$ to the $(i+1)^{th}$ component of this direct sum decomposition.  Then $$\eta(\xi \circ \Omega(\psi^i)) = G(\xi) \circ \eta(\Omega(\psi^i)) = G(\xi) \circ \begin{pmatrix} 0 \\ \vdots \\ 1 \\ \vdots \\ 0 \end{pmatrix}.$$  Thus with respect to the chosen decompositions, $G(\xi) : \Omega k_R^{n+1} \rightarrow k_R^{n+1}$ corresponds to the matrix with $\eta(\xi \circ \Omega(\psi^i))$ in the $(i+1)^{th}$ column.  Since $\psi^i : x^{m-1}A \rightarrow x^{m-1}A$ and $\Omega(\psi^i) : xA \rightarrow xA$ are both induced by left multiplication by $y^i$, while $\xi : xA \rightarrow x^{m-1}A$ is induced by multiplication by $x^{m-2}$, if $x^{m-2}$ commutes with $y$ we obtain $\xi \circ \Omega(\psi^i) = \psi^i \xi$ for each $i$.  Hence the $(i+1)^{th}$ column of $G(\xi)$ coincides with $$\eta(\psi^i \xi) = \begin{pmatrix} 0 \\ \vdots \\ f  \\ \vdots \\ 0 \end{pmatrix}$$with $f$ in the $(i+1)^{th}$ row as computed above.

  Alternatively, if $n=1$ we can write $\xi \Omega(\psi) = a \xi + b \psi \xi \in \stHom_A(\Omega T,T)$ for $a, b \in k$ with $b\neq 0$.  If $a \neq 0$ we would obtain $\psi \xi = a^{-1} \psi \xi \Omega(\psi)$, and thus $\psi \xi = a^{-2} \psi\xi\Omega(\psi)^2 = 0$, which is a contradiction.  Thus $\xi \Omega(\psi)$ equals a scalar multiple of $\psi \xi$.  As above, it follows that the second column of $G(\xi)$ would be $\bnc{0}{bf}$, and scaling the second generator of $\Omega(T_R)$ by $b^{-1}$ allows us to remove the $b$. 
 $\Box$ \\

\section{Spherical stable twists}
\setcounter{equation}{0}

The stable equivalences we construct in this section appear quite similar to the $0$-spherical twists introduced by Seidel and Thomas \cite{SeTh}, and we will see later that they arise naturally for group algebras of certain $2$-groups.

\begin{defin}\label{def:SphericalStableTwist}  Let $A$ be an elementary local symmetric $k$-algebra, and let $R$ be a unital subalgebra of $A$ such that ${}_R A$ and $A_R$ are free.  Let $\mu : A \otimes_R A \rightarrow A$ be induced by multiplication and write $C_{\mu}$ for the cone of $\mu$ in the category $\stlrp(A,A)$.  We write $\sigma_R : \stmod A \rightarrow \stmod A$ for the endo-functor induced by $-\otimes_A C_{\mu}$.

  If $R=k[x] \subset A$ for an element $x \in A$ with $x^m=0$, and $k \otimes_R A_R \cong k^2 \oplus R^l$ for some integer $l \geq 0$, then we call $\sigma_R$ a {\bf spherical stable twist}.  
\end{defin}

\begin{rmk}\label{rmk:relative-syzygy}  \textup{Suppose $M$ is an indecomposable, non-projective $A$-module.  Then $\sigma_R(M)$ is defined by the triangle $M \otimes_R A \rightarrow M \rightarrow \sigma_R(M) \rightarrow$, and it follows that we have a short exact sequence $\ses{\Omega \sigma_R(M)}{M \otimes_R A}{M}{}{}$.  On one hand, the map $M \otimes_R A \rightarrow M$ can be viewed as a relative $R$-projective cover of $M$, making $\Omega \sigma_R(M)$ a relative syzygy of $M$.  In particular, $\sigma_R(k)$ then resembles Alperin's construction of endo-trivial modules for certain $2$-groups \cite{Alp}.   Alternatively, when $m=2$ we have $M_R \in \add(k_R \oplus R)$, so that the map $M \otimes_R A \rightarrow M$ is a right $\add(T \oplus A)$-approximation of $M$.  Hence it is often straightforward to compute $\Omega \sigma_R(M)$ as the kernel of such an approximation. }
\end{rmk}

\begin{rmk}\label{rmk:spherical-objects}  \textup{ As noted in Remark~\ref{rmk:2dim}, the condition $k \otimes_R A_R \cong k^2 \oplus R^l$ is equivalent to $\stEndo_A(T) \cong k[\psi]/(\psi^2)$, where $T = k \otimes_R A$.   Thus $T$ resembles a {\it $0$-spherical object} as defined by Seidel and Thomas \cite{SeTh}, and we regard the functor $\sigma_R$ as an analogue of a $0$-sperical twist in this case.  In general, Seidel and Thomas define $n$-spherical objects in a bounded derived category $\K$, but their definition does not clearly carry over to other triangulated categories.  In particular, the finiteness condition 
\begin{quote}
(K2): for any $F \in \K$ both $\coprod_{i \in \mathbb{Z}} \K(E,F[i])$ and $\coprod_{i \in \mathbb{Z}} \K(F,E[i])$ are finite-dimensional,
\end{quote}
is not satisfied by any nonzero object $E$ in the stable module category of a self-injective algebra $\Lambda$, since we have $\Ext^i_{\Lambda}(E,\Lambda/J) \cong \stHom_{\Lambda}(\Omega^i E, \Lambda/J) \neq 0$ for all $i > 0$.  Moreover, this condition is essential to their definition of {\it twist functors} $T_E$, of which spherical twists are a special case.   As a substitute for this condition in $\stmod \Lambda$, we might require that $E$ is $\Omega$-periodic.  Then, if $E \cong E[n]$ say, we could replace $\coprod_{i \in \mathbb{Z}} \stHom(E,F[i])$ with $\coprod_{0 \leq i < n} \stHom(E,F[i])$, which is finite-dimensional and would still contain all the information about extensions involving $E$.  In fact, we have $T \cong T[2]$ in our situation, and we regard this isomorphism as an analogue of the condition (K2).}
\end{rmk}

\begin{therm}\label{theorem:SphericalStableTwist} A spherical stable twist $\sigma_R$ is an exact auto-equivalence of $\stmod A$.
\end{therm}

\noindent
{\it Proof.}  Since the left and right adjoints of $\sigma_R$ coincide, it suffices to show that $\sigma_R$ is fully faithful by Proposition 2.3.  Observe that $T[2] \cong T$ and $\C = \{T\} \cup T^{\perp}$ is a spanning class for $\stmod A$.  Furthermore $T^{\perp} = \{X \in \stmod A\ |\ \stHom_A(T,X) = 0\} = \{X \in \stmod A\ |\ \stHom_R(k,X_R) = 0\} = \{ X \in \stmod A\ |\ X_R \mbox{\ is\ projective}\}$.  Since $A_R$ is projective, $T^{\perp}$ is closed under (de-)suspensions, and it follows that there are no maps, in either direction, between $T[i]$ and $T^{\perp}$.  Thus, according to Propositon 2.2, to show that $\sigma_R$ is fully faithful, it suffices to check that it induces bijections $\stHom_A(T,T) \rightarrow \stHom_A(\sigma_R(T),\sigma_R(T))$, $\stHom_A(\Omega T,T) \rightarrow \stHom_A(\sigma_R(\Omega(T)),\sigma_R(T))$, and $\stHom_A(X,Y) \rightarrow \stHom_A(\sigma_R(X),\sigma_R(Y))$ for all $X, Y \in T^{\perp}$.

  If $X \in T^{\perp}$, then $\sigma_R(X)$ is defined as the cone of the natural map $X \otimes_R A \rightarrow X$.  Since $X_R$ is projective, so is $X \otimes_R A$ and thus we have an isomorphism $X \rightarrow \sigma_R(X)$ in $\stmod A$, and it is clearly natural with respect to $X \in T^{\perp}$.  Thus $\sigma_R$ is fully faithful on $T^{\perp}$.

  We now consider $\sigma_R(T)$.   For ease of notation, we let $F$ and $G$ be the functors induced by induction and restriction respectively between $\stmod R$ and $\stmod A$.  We let $\eta : \stHom_A(FX,Y) \stackrel{\cong}{\longrightarrow} \stHom_R(X,GY)$ denote the adjugant.  Observe that for any $X \in \stmod A$, the counit $\delta_X = \eta^{-1}(1_{GX}) : FGX \rightarrow X$ coincides with the natural map $X \otimes_R A \rightarrow X$ induced by tensoring the map $\mu$ with $X_A$ on the left.   By Lemma 4.1 $GT \cong k^2$ and with respect to the basis given in the proof, we have $\eta(1) : k \rightarrow GT$ corresponding to the map $\bnc{1}{0}$ and $\eta(\psi) : k \rightarrow GT$ corresponding to $\bnc{0}{1}$.   Hence $F(\eta(1)) = \bnc{1_T}{0} : T \rightarrow FGT \cong T^2$ and $F(\eta(\psi)) = \bnc{0}{1_T} : T \rightarrow FGT \cong T^2$.  Now $\delta_T : FGT \cong T^2 \rightarrow T$ has the form $(u\ v)$ for $u, v \in \stEndo_A(T)$.  Since $1_T = \delta_T F(\eta(1_T)) = \delta_T \bnc{1_T}{0} = u$ and $\psi = \delta_T F(\eta(\psi)) = \delta_T \bnc{0}{1_T} =v$, we see that $\delta_T = (1_T, \psi)$ with respect to this decomposition $FGT \cong T^2$.  Now consider the induced isomorphism of split triangles
  \begin{equation} \vcenter{ \xymatrix{ \sigma_R(T)[-1] \ar[r] \ar[d]^{\cong} & FG(T) \ar[r] \ar[d]^{\cong} & T \ar[r] \ar@{=}[d] & \\ T \ar[r]_{\bnc{-\psi}{1}} & T^2 \ar[r]_{(1\ \psi)} & T \ar[r] & }} 
  \end{equation}
   which shows $\sigma_R(T)[-1] \cong T$ and thus $\sigma_R(T) \cong \Sigma(T)$.  By Lemma 4.1 we have $FG(\psi) = \bnc{0 \ 0}{1_T \ 0} : FGT \rightarrow FGT$, and using the above isomorphism of triangles we compute $\sigma_R(\psi) = \Sigma(-\psi)$:
    
    $$\xymatrixcolsep{3.5pc} \xymatrix{ T \ar[r]^{\bnc{-\psi}{1}} \ar@{-->}[d]^{-\psi} & T^2 \ar[r]^{(1\ \psi)} \ar[d]^{\left( \begin{smallmatrix} 0 & 0 \\ 1_T & 0 \end{smallmatrix} \right)} &  T \ar[r]^(.35){0} \ar[d]^\psi & \sigma_R(T) \cong \Sigma(T) \ar@{-->}[d]^{\Sigma(-\psi)} \\ T \ar[r]_{\bnc{-\psi}{1}} & T^2 \ar[r]_{(1\ \psi)} &  T \ar[r]^(.35){0} & \sigma_R(T) \cong \Sigma(T)}.$$ 
    Thus $\sigma_R$ induces an isomorphism $\stEndo_A(T) \cong \stEndo_A(\sigma_R(T))$ as required.   

Next we consider the action of $\sigma_R$ on the map $\xi \in \stHom_A(\Omega T,T)$ as in Lemma~\ref{lemma:RestrictionXi}.  As in the proof of that lemma, we have a nonzero scalar $b \in k$ such that $\xi \Omega(\psi) = b \psi \xi$.  If we suspend the isomorphism of triangles (5.1), we can use it to compute $\sigma_R(\xi) = \Sigma(b\xi)$ as in the diagram below, where the map $\Omega T^2 \rightarrow T^2$ is $FG(\xi)$ according to Lemma~\ref{lemma:RestrictionXi}.
$$\xymatrixcolsep{3.5pc} \xymatrix{ \Omega T \ar[r]^{\bnc{-\Omega(\psi)}{1}} \ar@{-->}[d]^{b\xi} & \Omega T^2 \ar[r]^{(1\ \Omega(\psi))} \ar[d]^{\left( \begin{smallmatrix} \xi & 0 \\ 0 & b\xi \end{smallmatrix} \right)} &  \Omega T \ar[r]^(.35){0} \ar[d]^\xi & \sigma_R(\Omega(T)) \cong \Sigma(\Omega T) \ar@{-->}[d]^{\Sigma(b\xi)} \\ T \ar[r]_{\bnc{-\psi}{1}} & T^2 \ar[r]_{(1\ \psi)}  & T \ar[r]^(.35){0} & \sigma_R(T) \cong \Sigma(T)}.$$ 
It follows that $\sigma_R$ induces an isomorphism from $\stHom_A(\Omega T,T)$ to $\stHom_A(\sigma_R(\Omega T), \sigma_R(T))$ as required.  By Proposition 2.2, we conclude that $\sigma_R$ is an equivalence.  $\Box$ \\

It is not difficult to see that a spherical stable twist $\sigma_R$ is (usually) not isomorphic to any power of $\Omega$ on $\stmod A$, simply by comparing the dimensions of $\sigma_R(k_A)$ and $\Omega^n(k_A)$ for $n \in \mathbb{Z}$.  

\begin{propos} Let $d = \dim_k A$, and assume that $d/m = \rk_R A > 2$.  Then $\sigma_R(k_A) \ncong \Omega^n(k_A)$ for all $n \in \mathbb{Z}$.  Hence $\sigma_R$ is not induced by an auto-equivalence of $D^b(\rmod A)$.
\end{propos}

\noindent
{\it Proof.}  We consider the stable Grothendieck group $\ul{K}_0(A) \cong \mathbb{Z}/d\mathbb{Z}$.  Any auto-equivalence of $\stmod A$ induces an automorphism of $\ul{K}_0(A)$, and the syzygy functor $\Omega$ corresponds to multiplication by $-1$ on $\mathbb{Z}/d\mathbb{Z}$.  Since $\sigma_R(k_A)$ is defined to be the cosyzygy of the kernel of the natural map $k \otimes_R A \rightarrow k$, we see that the image of $\sigma_R(k_A)$ in the stable Grothendieck group is $1 - \frac{d}{m} \neq \pm 1$. The first statement is now immediate, and the second follows from Corollary~\ref{coro:no-lift}.  $\Box$ \\

From the above proof, we see that $\sigma_R$ corresponds to multiplication by $1 - d/m$ in $\ul{K}_0(A) \cong \mathbb{Z}/d\mathbb{Z}$.  This yields some information about powers of $\sigma_R$, but it is generally a more difficult problem to determine whether or not any higher powers of $\sigma_R$ coincide with powers of $\Omega$.  Two different possibilities are illustrated by the dihedral and semidihedral algebras that we examine in the last section.

\section{$\mathbb{P}^n$ stable twists}
\setcounter{equation}{0}

We now define additional auto-equivalences of $\stmod A$ that are analogous to the $\mathbb{P}^n$ twists introduced by Huybrechts and Thomas \cite{HuTh}.  We continue the notation introduced in the last section.  Namely $A$ is an elementary, local symmetric $k$-algebra and $R$ is the unital subalgebra generated by an element $x \in \rad A$ with $x^m=0$ and $x^{m-1} \neq 0$ for some $m \geq 2$.  We also assume that ${}_R A$ and $A_R$ are both free, and we set $T = k \otimes_R A$.

\begin{defin} Let $R$ and $A$ be as above and suppose $y \in A$ commutes with $x$.  Let $H : A\otimes_R A \rightarrow A \otimes_R A$ be the $(A,A)$-bimodule map sending $1 \otimes 1$ to $y\otimes 1 - 1\otimes y$ and $\mu : A \otimes_R A \rightarrow A$ the map induced by multiplication.  Let $C_H$ be the cone of $H$ in $\stlrp(A,A)$, and observe that $\mu$ factors through $C_H$ since $\mu H = 0$.  Let $Q$ be the cone of the induced map from $C_H$ to $A$ in $\stlrp(A,A)$.  
\begin{equation}\vcenter{ \xymatrix{ 
A \otimes_R A \ar[r]^{H} \ar@{-->}[dr]^0 & A \otimes_R A \ar[r] \ar[d]^{\mu} &  C_H \ar[d] \ar[r] & (A\otimes_R A)[1] \\ & A \ar@{=}[r] \ar[d] & A \ar[d] \\ (A \otimes_R A)[1] \ar[r] & C_{\mu} \ar[r] & Q \ar[r] & (A \otimes_R A)[2]} }\end{equation}

\noindent
We write $\rho_{R,y} : \stmod A \rightarrow \stmod A$ for the endo-functor induced by $-\otimes_A Q$.  If $\stEndo_A(A/xA) \cong k[\psi]/(\psi^{n+1})$ with $\psi$ induced by left multiplication by $y$, then we call $\rho_{R,y}$ a {\bf  $\mathbb{P}^n$ stable twist}.  
\end{defin}

\begin{therm} A $\mathbb{P}^n$ stable twist $\rho=\rho_{R,y}$ is an auto-equivalence of $\stmod A$.
\end{therm}

The proof is similar to that of Theorem 5.2.  In particular, we start by analyzing the action of $\rho$ on the spanning class $\C = \{T\}\cup T^{\perp}$.  In Lemmas 6.3 and 6.4, we show that $\rho : \stHom_A(X,Y[i]) \rightarrow \stHom_A(\rho(X),\rho(Y[i]))$ is bijective for any $X, Y \in \C$.  The theorem then follows easily by Propositions 2.2 and 2.3.

\begin{lemma} Restricted to $T^{\perp}$, $\rho$ is isomorphic to the identity functor.
\end{lemma}

\noindent
{\it Proof.}  As in the proof of Theorem 5.2, $T^{\perp}$ coincides with the full subcategory of $\stmod A$ consisting of $A$-modules that are free over $R$.  By construction, for $X \in T^{\perp}$, $X \otimes_A C_H$ is projective.  Hence the map $A \rightarrow Q$ induces isomorphisms $X \rightarrow X \otimes_A Q$ in $\stmod A$. $\Box$\\

\begin{lemma} Restricted to $\add(T \oplus T[-1])$, $\rho$ is isomorphic to $\Sigma^2$.
\end{lemma}

\noindent
{\it Proof.}  Tensoring the triangle $A\otimes_R A \stackrel{H}{\longrightarrow} A\otimes_R A \longrightarrow C_H \rightarrow$ with $T_A$ and using Lemma 4.1, we obtain the top triangle in the morphism of triangles below.
\begin{equation}\vcenter{\xymatrixrowsep{6.0pc} \xymatrixcolsep{7.0pc} \xymatrix{ T^{n+1} \ar[d]_{\left( \begin{smallmatrix} 1 & & & -y^{n} \\ & \ddots & & \vdots \\ & & 1 & -y \\ & & & 1 \end{smallmatrix} \right)} \ar[r]^{\left(\begin{smallmatrix} -y \\ 1 & -y \\ & \ddots & \ddots \\ & & 1 & -y \end{smallmatrix} \right)} & T^{n+1} \ar[r] \ar[d]^{\left(\begin{smallmatrix} 0 & 1 & y & \cdots & y^{n-1}\\ & & \ddots & \ddots & \vdots \\ & & & 1 & y \\ & & & & 1 \\ 1 & y & y^2 & \cdots & y^{n} \end{smallmatrix} \right)} & T \otimes_A C_H \ar@{-->}[d]^{\cong} \ar[r] & \\ T^{n+1} \ar[r]_{\left( \begin{smallmatrix} I_{n} & 0 \\ 0 & 0 \end{smallmatrix} \right)} & T^{n+1} \ar[r]_{\left(\begin{smallmatrix} 0 & \cdots & 1 \\ 0 & \cdots & 0 \end{smallmatrix}\right)} & T \oplus T[1] \ar[r] & }}\end{equation}
As the first two vertical maps are isomorphisms, so is the third.  Furthermore, as in the proof of Theorem 5.2, the map $\mu : A \otimes_R A \rightarrow A$ induces the map $1_T \otimes \mu : (1_T\ \psi \ \cdots \ \psi^{n}): T^{n+1} \rightarrow T$.  Under the above isomorphism of triangles this map corresponds to the projection of $T^{n+1}$ onto its $(n+1)^{th}$ factor.  Consequently, the induced map $T \oplus T[1] \cong T \otimes_A C_H \rightarrow T \otimes_A A \cong T$ is the projection onto $T$.  As $\rho(T)$ is defined as the cone of this map, we see that $\rho(T) \cong T[2]$:

\begin{equation} \vcenter{\xymatrix{ \rho(T)[-1] \ar[r] \ar@{-->}[d]^{\cong} & T \otimes_A C_H \ar[r] \ar[d]^{\cong} & T \ar[r] \ar@{=}[d] & \rho(T) \ar@{-->}[d]^{\cong} \\ T[1] \ar[r]^(.45){\bnc{0}{1}} & T \oplus T[1]  \ar[r]^(.65){(1\ 0)} & T \ar[r]^0 & T[2]}} \end{equation}

  Of course it follows now that $\rho(T[-1]) \cong T[1]$, but notice that we can also suspend the above argument.  We point out that $\Sigma(\psi)$ is also induced by multiplication by $y$ on $T[-1] \cong A/x^{m-1}A \cong xA$.

We now compute the effect of $\rho$ on $\stHom_A(T,T) = k[\psi]/(\psi^{n+1})$ and on $\stHom_A(\Omega T,T)$.  When we tensor the map $\psi : T \rightarrow T$ with the triangle $A \otimes_R A \rightarrow A \otimes_R A \rightarrow C_H$ and pass to the isomorphic triangle in (6.2), we see that the resulting endomorphism of the left-most $T^{n+1}$ is $$\begin{pmatrix}1 & & & -y^{n} \\ & \ddots & & \vdots \\ & & 1 & -y \\ & & & 1 \end{pmatrix} \begin{pmatrix} 0 \\ 1 & 0 \\ & \ddots \\ & & 1 & 0 \end{pmatrix} \begin{pmatrix} 1 & & & -y^{n} \\ & \ddots & & \vdots \\ & & 1 & -y \\ & & & 1 \end{pmatrix}^{-1} = \begin{pmatrix} 0 & &\cdots & -y^{n} & 0 \\ 1 & 0 & \cdots & -y^{n-1} & 0 \\ &  \ddots & & \vdots & \vdots \\ & & 1 & -y & 0 \\ 0 & \cdots & 0 & 1 & y\end{pmatrix}.$$
Restricted to the $(n+1)^{th}$ factor of $T$ we thus obtain the endomorphism $\psi$, and hence we have $\Sigma(\psi)$ for the $T[1]$-component of the induced endomorphism of $T \otimes_A C_H \cong T \oplus T[1]$.  Finally, this yields $\Sigma^2(\psi)$ for the corresponding endomorphism of $\rho(T) \cong T[2]$.

We now compute the effect of $\rho$ on the map $\xi \in \stHom_A(\Omega T, T)$.  By Lemma~\ref{lemma:RestrictionXi}, $FG(\xi) : (\Omega T)^{n+1} \rightarrow T^{n+1}$ will be a given by a diagonal matrix with $\xi$ along the diagonal.  This representation of $FG(\xi)$ does not change when we pass to the isomorphic triangles as in (6.2), and hence we see that the induced map $\Omega(T) \oplus \Omega(T)[1] \rightarrow T \oplus T[1]$ is $\left( \begin{smallmatrix} \xi & 0 \\ 0 & \Sigma{\xi} \end{smallmatrix} \right)$.  For the second cone construction, we use triangles as in (6.3) and its de-suspension to compute $\rho(\xi) = \Sigma^2(\xi)$:

\begin{equation} \vcenter{\xymatrix{ \Omega(T)[1] \ar[r]^(.45){\bnc{0}{1}} \ar@{-->}[d]^{\Sigma \xi} & \Omega T \oplus \Omega T [1] \ar[r]^(.65){(1\ 0)} \ar[d]^{\left( \begin{smallmatrix} \xi & 0 \\ 0 & \Sigma{\xi} \end{smallmatrix} \right)} & \Omega T \ar[r]^(.4){0} \ar[d]^{\xi} & \Omega T [2] \cong \rho(\Omega T) \ar@{-->}[d]^{\Sigma^2(\xi)} \\ T[1] \ar[r]^(.45){\bnc{0}{1}} & T \oplus T[1]  \ar[r]^(.65){(1\ 0)} & T \ar[r]^(.4){0} & T[2] \cong \rho(T)}} \end{equation}

Finally, since all maps between $T$ and its suspensions are generated by $\psi$ and $\xi$, on which $\rho$ acts as $\Sigma^2$, we conclude that the restriction of $\rho$ to $\add(T \oplus T[-1])$ is isomorphic to $\Sigma^2$.  $\Box$\\

\begin{propos} Let $\rho = \rho_{R,y}$ be a $\mathbb{P}^n$ stable twist.  Then $$\rho_{R,y}(k) \cong \Omega^{-2}(\Omega^2_R(k)),$$ where $\Omega_R (M)$  denotes the kernel of a relatively $R$-projective cover of $M$.  
\end{propos}

\noindent
{\it Proof.}  Tensoring the triangles in the diagram (6.1) with $k$ yields the following commutative diagram in which the rows and columns are triangles

\begin{equation}\vcenter{\xymatrix{ T \ar@{=}[d] \ar[r] & \Omega_R(k) \ar[d] \ar[r] & \rho_R(k)[-1] \ar[r] \ar[d] & T[1] \ar@{=}[d] \\ T \ar[r]^{\psi} \ar@{-->}[dr]^0 & T \ar[r] \ar[d] & k \otimes_A C_H \ar[d] \ar[r] & T[1] \\ & k \ar@{=}[r] \ar[d] & k \ar[d] \\ T[1] \ar[r] & \Omega_R(k)[1] \ar[r] & \rho_R(k) \ar[r] & T[2]}} \end{equation}
We claim that the map $T \rightarrow \Omega_R(k)$ in the top left is a relatively $R$-projective cover.  In fact, when we restrict this map to $\stmod R$, we get a split epimorphism $k^{n+1} \rightarrow k^n$ by Lemma 4.1.  It follows that for any $R$-module $M$, the induced map $\stHom_R(M, T_R) \rightarrow \stHom_R(M, (\Omega_R(k))_R)$ is onto, and hence $\stHom_A(FM, T) \rightarrow \stHom_A(FM,\Omega_R(k))$ is also onto; i.e., the map $T \rightarrow \Omega_R(k)$ is a right $F(\stmod R)$-approximation.  It follows that $\Omega^2_R(k) \cong \rho_R(k)[-2]$  as desired. $\Box$ \\

In the context of the derived category of coherent sheaves of a smooth projective variety, Huybrechts and Thomas have shown that a $\mathbb{P}^1$-twist coincides with the square of a spherical twist \cite{HuTh}.  The analagous statement holds in our context, at least up to a Morita auto-equivalence, whenever $m=2$ and a $\mathbb{P}^1$-twist is defined.

\begin{propos}  Let $A$ be a split, local, symmetric $k$-algebra which is free on either side over the subalgebra $R=k[x]$ where $x^2=0$.  Suppose that $\stEndo_A(k \otimes_R A) \cong k[\psi]/(\psi^2)$ with $\psi$ corresponding to left multiplication by some $y \in A$ that commutes with $x$.  Then $\rho_{R,y}(k) \cong \sigma_R^2(k)$.  In particular, there is a Morita auto-equivalence $F$ of $\mod A$ (i.e., $F = -\otimes_A {}_1 A_{\sigma}$ for a twisted bimodule ${}_1 A_{\sigma}$ with $\sigma \in \mbox{Out}(A)$) such that $\rho_{R,y} \cong \sigma_R^2 \circ F$ as functors on $\stmod A$.

\end{propos}

\noindent
{\it Proof.}  Notice we may replace the $\Omega_R(k)[1]$ entry in diagram (6.5) with $\sigma_R(k)$.  Now $\stHom_A(T[1], \sigma_R(k)) \cong \stHom_A(\sigma_R(T), \sigma_R(k)) \cong \stHom_A(T,k)$ is one-dimensional and $T[1] \cong T$.  Thus the map $T[1] \rightarrow \sigma_R(k)$ in the bottom row of the diagram is either $0$ or else a right $\add(T)$-approximation.  It must be nonzero since $\Sigma \psi$ factors through it.  Thus we have $\rho_R(k) \cong \sigma_R^2(k)$.  The final statement now follows from Proposition~\ref{prop:LinckelmannTheorem}. $\Box$\\

\section{Examples and questions}
\setcounter{equation}{0}

\subsection{Dihedral groups and algebras}

Let $k$ be any field and set $A = k\gen{x,y}/(x^2,y^2,(xy)^q-(yx)^q)$ for some $q \geq 2$.  In general, $A$ is a local, symmetric, special biserial $k$-algebra of dimension $4q$ with $k$-basis $$\{(xy)^i, (xy)^ix, (yx)^{i+1}, (yx)^iy\ |\ 0 \leq i < q\}.$$
For arbitrary fields $k$ and integers $q \geq 2$, the indecomposable $A$-modules have been completely described by Ringel \cite{Rin} and the Auslander-Reiten quiver can be worked out as in \cite{Erd}, Lemma II.7.6 (see also \cite{RAC1}, Sections 4.11 and 4.17).    If $\chr(k) = 2$ and $q$ is a power of $2$, then $A$ is isomorphic to the group algebra over $k$ of the dihedral group of order $4q$.

We set $R = k[x] \subset A$, and it is clear that ${}_R A$ and $A_R$ are free $R$-modules of rank $2q$ with bases given by $\{(yx)^iy^j\ |\ 0 \leq i < q, j=0,1\}$ and $\{y^j(xy)^i\ |\ 0 \leq i < q, j=0,1\}$, respectively.  Furthermore, as an $(R,R)$-bimodule, $A$ decomposes as two copies of the regular bimodule $R$, generated by $1$ and $(yx)^{q-1}y$, and $q-1$ copies of the projective bimodule $R \otimes_k R$, generated by $(yx)^iy$ for $0 \leq i < q-1$; i.e., ${}_R A_R \cong R^2 \oplus (R \otimes_k R)^{q-1}$.  Following the notation of Section 4, we see that $T = k \otimes_R A \cong A/xA$ is a uniserial module of dimension $2q$, corresponding to the word $(yx)^{q-1}y$.  Clearly, $T_R = k \otimes_R A_R \cong k \otimes_R (R^2 \oplus (R \otimes_k R)^{q-1}) \cong k^2 \oplus R^{q-1}$.  Thus the hypotheses of Definition 5.1 are satisfied, and we see that the spherical twist $\sigma_R$ gives an auto-equivalence of $\stmod A$.  Applying $\sigma_R$ to the simple $A$-module $k$, we see that $\sigma_R(k)$ and $\sigma_R^2(k)$ are the string modules corresponding to the words $x(xy)^{1-q}y^{-1}$ and $x(xy)^{1-q}y^{-1}x(xy)^{1-q}y^{-1} $, respectively, which have graphs
$$\begin{array}{ccc} 
\xymatrixrowsep{0.75pc} \xymatrixcolsep{0.75pc} 
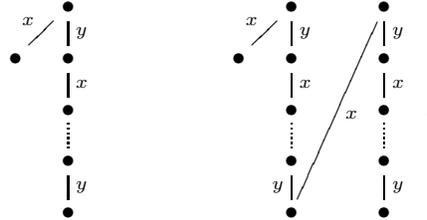
\xymatrix{ & \bullet \ar@{-}[dl]_x \ar@{-}[d]^y \\ \bullet &  \bullet \ar@{-}[d]^x \\ &  \bullet \ar@{.}[d] \\  & \bullet \ar@{-}[d]^y \\  & \bullet} & \hspace{1cm} & 
\xymatrixrowsep{0.75pc} \xymatrixcolsep{0.75pc} \xymatrix{ & \bullet \ar@{-}[dl]_x \ar@{-}[d]^y & & \bullet \ar@{-}[ddddll]^x \ar@{-}[d]^y \\
 \bullet  & \bullet \ar@{-}[d]^x & & \bullet \ar@{-}[d]^x \\
   & \bullet \ar@{.}[d] & & \bullet \ar@{.}[d]\\
   & \bullet \ar@{-}[d]_y & & \bullet \ar@{-}[d]^y   \\ 
   &  \bullet  & & \bullet } 
 \end{array}.$$

When $A = kD_{4q}$, $\sigma_R(k)$ coincides with $\Omega^{-1}(L)$ for the endo-trivial module $L$ described in \cite{CarThe1}, Section 5.  It is shown there that the group $T(D_{4q})$ of endo-trivial modules is isomorphic to $\mathbb{Z}^2$, generated by the classes of $\Omega(k)$ and $L$.  In this case, the fact that no power of $L$ coincides with a power of $\Omega(k)$--or, equivalently, that no power of $\sigma_R$ coincides with a power of $\Omega$--can be seen by restricting to two different elementary abelian $2$-subgroups of rank $2$.  In the stable category of one such subgroup, $\sigma_R(k)$ restricts to $k$, while in the other it restricts to $\Omega^{-2}(k)$.  

Similar arguments do not appear to work in general for showing that no nonzero power of $\sigma_R$ can coincide with a power of $\Omega$ on $\stmod A$.  However, we can see that this is indeed still true by considering the action of $\sigma_R$ on the stable AR-quiver of $A$.  We focus on the two components of the stable AR-quiver $\mathcal{C}_k$ and $\mathcal{C}_A$ containing the simple $A$-module and $\rad A$, respectively.  Each of these components has the shape $\mathbb{Z}A_{\infty}^{\infty}$, and $\Omega$ induces isomorphisms of translation quivers $\mathcal{C}_k \rightarrow \mathcal{C}_A$ and $\mathcal{C}_A \rightarrow \mathcal{C}_k$ with $\Omega^2 = \tau$ since $A$ is symmetric.  It can be seen that $\sigma_R(k)$ is one of the two immediate successors of $k$ in $\mathcal{C}_k$  (the other successor is in fact $\sigma_{R'}(k)$ for $R'=k[y] \subseteq A$).  Furthermore, one checks that $\sigma_R^2(k) \ncong \Omega^{-2}(k)=\tau^{-1}(k)$.  Since the stable equivalence $\sigma_R$ commutes with $\Omega$ on objects, it follows that $\sigma_R$ induces an automorphism of $\mathcal{C}_k$ that can be described as shifting down and to the right one arrow.  Again since $\sigma_R$ commutes with $\Omega$, $\sigma_R$ induces an automorphism of $\mathcal{C}_A$ of the same form.  We sketch a portion of the component $\C_k$ below.

In particular, we see that no nonzero power of $\sigma_R$ can be isomorphic to a power of $\Omega$.  Moreover, in the case where $A = kD_{4q}$ (with $q$ a power of $2$), we see that these two components of the stable AR-quiver consist precisely of all the endo-trivial $A$-modules.  This fact has been previously observed by Bessenrodt \cite{Bes}.  While it appears interesting, it is probably not indicative of what happens more generally.

$$\xymatrixrowsep{2.0pc} \xymatrixcolsep{4.0pc} 
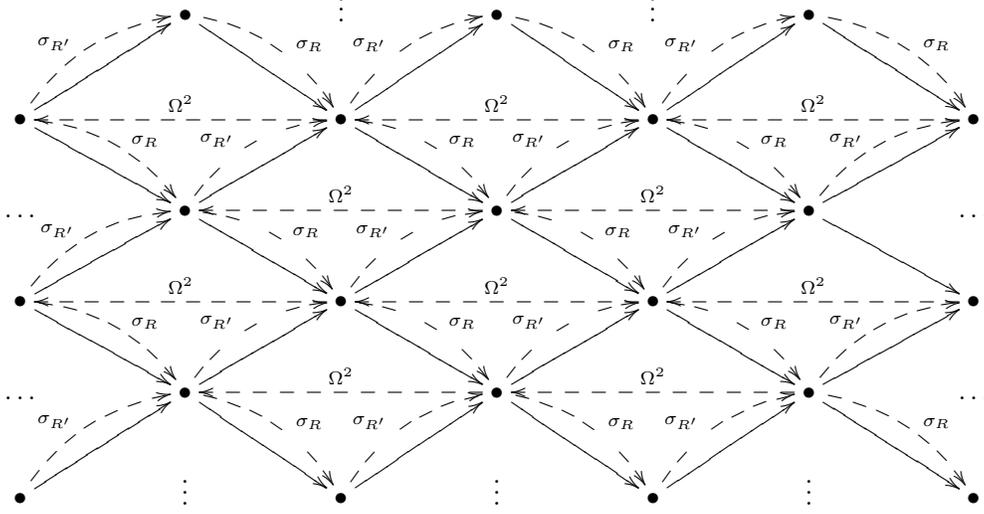
\xymatrix{ & \bullet \ar[dr] \ar@{-->}@/^.75pc/[dr]^(.6){\sigma_R} & \vdots & \bullet \ar[dr] \ar@{-->}@/^.75pc/[dr]^(.6){\sigma_R} & \vdots & \bullet \ar[dr] \ar@{-->}@/^.75pc/[dr]^(.6){\sigma_R}& \\
 \bullet   \ar[ur] \ar@{-->}@/^.75pc/[ur]^(.4){\sigma_{R'}} \ar[dr] \ar@{-->}@/^.75pc/[dr]^(.6){\sigma_R} & & \bullet \ar[ur] \ar@{-->}@/^.75pc/[ur]^(.4){\sigma_{R'}} \ar[dr] \ar@{-->}@/^.75pc/[dr]^(.6){\sigma_R} \ar@{-->}[ll]_{\Omega^2}  & & \bullet \ar[ur] \ar@{-->}@/^.75pc/[ur]^(.4){\sigma_{R'}} \ar[dr] \ar@{-->}@/^.75pc/[dr]^(.6){\sigma_R} \ar@{-->}[ll]_{\Omega^2} & & \bullet \ar@{-->}[ll]_{\Omega^2} \\
  \hdots & \bullet \ar[ur] \ar@{-->}@/^.75pc/[ur]^(.4){\sigma_{R'}}  \ar[dr] \ar@{-->}@/^.75pc/[dr]^(.6){\sigma_R} & & \bullet \ar[ur] \ar@{-->}@/^.75pc/[ur]^(.4){\sigma_{R'}}  \ar[dr] \ar@{-->}@/^.75pc/[dr]^(.6){\sigma_R} \ar@{-->}[ll]_{\Omega^2} & & \bullet \ar[ur] \ar@{-->}@/^.75pc/[ur]^(.4){\sigma_{R'}} \ar[dr] \ar@{-->}[ll]_{\Omega^2} & \hdots \\ 
  \bullet  \ar[ur] \ar@{-->}@/^.75pc/[ur]^(.4){\sigma_{R'}} \ar[dr]  \ar@{-->}@/^.75pc/[dr]^(.6){\sigma_R} & & \bullet \ar[ur] \ar@{-->}@/^.75pc/[ur]^(.4){\sigma_{R'}} \ar[dr] \ar@{-->}@/^.75pc/[dr]^(.6){\sigma_R}   \ar@{-->}[ll]_{\Omega^2} & & \bullet \ar[ur] \ar@{-->}@/^.75pc/[ur]^(.4){\sigma_{R'}} \ar[dr] \ar@{-->}@/^.75pc/[dr]^(.6){\sigma_R} \ar@{-->}[ll]_{\Omega^2} & & \bullet \ar@{-->}[ll]_{\Omega^2} \\ 
  \hdots & \bullet \ar[ur] \ar@{-->}@/^.75pc/[ur]^(.4){\sigma_{R'}} \ar[dr] \ar@{-->}@/^.75pc/[dr]^(.6){\sigma_R} & & \bullet \ar[ur] \ar@{-->}@/^.75pc/[ur]^(.4){\sigma_{R'}} \ar[dr] \ar@{-->}@/^.75pc/[dr]^(.6){\sigma_R} \ar@{-->}[ll]_{\Omega^2} & & \bullet \ar[ur] \ar@{-->}@/^.75pc/[ur]^(.4){\sigma_{R'}} \ar[dr] \ar@{-->}@/^.75pc/[dr]^(.6){\sigma_R} \ar@{-->}[ll]_{\Omega^2} & \hdots \\ 
  \bullet \ar[ur]  \ar@{-->}@/^.75pc/[ur]^(.4){\sigma_{R'}} & \vdots & \bullet \ar[ur] \ar@{-->}@/^.75pc/[ur]^(.4){\sigma_{R'}} & \vdots & \bullet \ar[ur] \ar@{-->}@/^.75pc/[ur]^(.4){\sigma_{R'}} & \vdots & \bullet }$$

In a forthcoming article \cite{StPic}, we will further apply this analysis to describe the stable Picard groups, consisting of isomorphism classes of auto-equivalences of $\stmod A$, for local dihedral algebras $A$.

\subsection{Semidihedral groups and algebras} Let $k$ be any field and set $$A = k\gen{x,y}/(x^2,y^4,y^2-(xy)^{q-1}x-\delta(yx)^q,(xy)^q-(yx)^q)$$ for some $q \geq 2$ and $\delta \in k$.  In general, $A$ is a local, symmetric $k$-algebra of dimension $4q$ with $k$-basis $$\{(xy)^i, (xy)^ix, (yx)^{i+1}, (yx)^iy\ |\ 0 \leq i < q\}.$$
 If $\chr(k) = 2$, $q = 2^n$ for $n\geq 2$ and $\delta=1$, then $A$ is isomorphic to the group algebra over $k$ of the semidihedral group of order $4q$.

As for the dihedral algebras above, we may set $R = k[x] \subset A$, and it is clear that ${}_R A$ and $A_R$ are free $R$-modules of rank $2q$ with bases given by $\{(yx)^iy^j\ |\ 0 \leq i < q, j=0,1\}$ and $\{y^j(xy)^i\ |\ 0 \leq i < q, j=0,1\}$, respectively.  Furthermore, as an $(R,R)$-bimodule, $A$ decomposes as two copies of the regular bimodule $R$, generated by $1$ and $(yx)^{q-1}y$, and $q-1$ copies of the projective bimodule $R \otimes_k R$, generated by $(yx)^iy$ for $0 \leq i < q-1$; i.e., ${}_R A_R \cong R^2 \oplus (R \otimes_k R)^{q-1}$.  Following the notation of Section 4, we see that $T = k \otimes_R A \cong A/xA$ is a uniserial module of dimension $2q$, corresponding to the word $(yx)^{q-1}y$.  Clearly, $T_R = k \otimes_R A_R \cong k \otimes_R (R^2 \oplus (R \otimes_k R)^{q-1}) \cong k^2 \oplus R^{q-1}$.  Thus the hypotheses of Definition 5.1 are satisfied, and we see that the spherical twist $\sigma_R$ gives an auto-equivalence of $\stmod A$.  Applying $\sigma_R$ to the simple $A$-module $k$, we see that $\sigma_R(k)$ is the string module corresponding to the word $x(xy)^{1-q}y^{-1}$ as in the dihedral case.

When $A = kSD_{4q}$, $\sigma_R(k)$ coincides with $\Omega^{-1}(L)$ for the endo-trivial module $L$ described in \cite{CarThe1}, Section 7.  It is shown there that the group $T(SD_{4q})$ of endo-trivial modules is isomorphic to $\mathbb{Z} \oplus \mathbb{Z}/2\mathbb{Z}$, generated by the classes of $\Omega(k)$ and $\Omega(L)$, which has order $2$ (see Theorem 7.1 of \cite{CarThe1}).   Hence, in this case we have $\sigma_R^2(k) \cong \Omega^{-4}(k)$ and it follows that $(-\otimes_k \sigma_R(k))^2 \cong \Omega^{-4}$ on $\stmod A$.

  It is not clear whether the analagous isomorphism $\sigma_R^2 \cong \Omega^{-4}$ holds for all local semidihedral algebras.  However, we can again consider the automorphism of the stable AR-quiver of $A$ induced by $\sigma_R$.  In II.10 of \cite{Erd}, Erdmann shows that the stable AR-quiver of $A$ consists of (infinitely many) $\mathbb{Z}A_{\infty}^{\infty}$ and $\mathbb{Z}D_{\infty}$ components, as well as tubes of rank $1$ and $2$.  The simple $A$-module $k$ lies on the mouth of a $\mathbb{Z}D_{\infty}$ component.  The almost split sequence starting in $k$ can be computed by applying $\Omega^{-1}$ to the almost-split sequences starting in $\rad A$, which has the form $\ses{\rad A}{A \oplus \rad A/ \soc A}{A/\soc A}{}{}$.  It follows that $\ses{k}{\Omega^{-1}(\rad A/\soc A)}{\Omega^{-2}(k)}{}{}$ is an almost-split sequence.  Observe that $\rad A/ \soc A$ is indecomposable.  There is another almost split sequence $$\ses{\Omega^2\sigma_R(k)}{\Omega^{-1}(\rad A/\soc A)}{\sigma_R(k)}{}{},$$ which is $\Omega^{-1}$ applied to the sequence in II.9.5 of \cite{Erd}.  Thus the component of $k$ is preserved by $\sigma_R$ which has the effect of shifting vertices one unit to the right, except for the two leaves of each sectional $D_{\infty}$ subtree, which are swapped and then shifted one unit to the right.  It is clear that $\sigma_R^2(X) \cong \Omega^{-4}(X)$ for each indecomposable $A$-module $X$ in this component, and for most such $X$ (those not on the mouth), we even have $\sigma_R(X) \cong \Omega^{-2}(X)$.

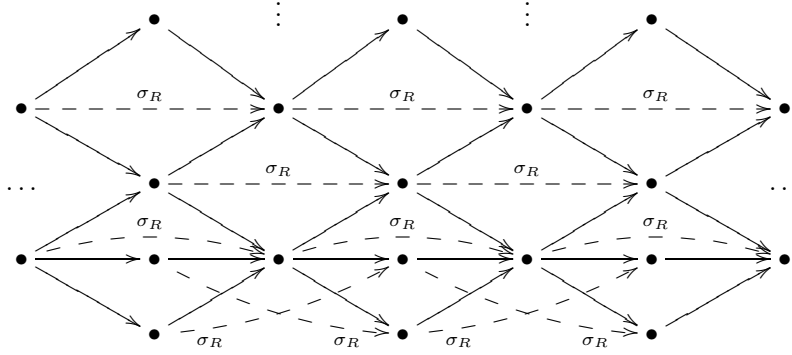
  $$ \xymatrixrowsep{1.5pc} \xymatrixcolsep{3.0pc} \xymatrix{ & \bullet \ar[dr] & \vdots & \bullet \ar[dr] & \vdots & \bullet \ar[dr] \\
  \bullet \ar[ur] \ar[dr] \ar@{-->}[rr]^{\sigma_R} & & \bullet \ar[ur] \ar[dr] \ar@{-->}[rr]^{\sigma_R} & & \bullet \ar[ur] \ar[dr] \ar@{-->}[rr]^{\sigma_R} & & \bullet \\
  \hdots &  \bullet \ar[ur] \ar[dr] \ar@{-->}[rr]^{\sigma_R} & & \bullet \ar[ur] \ar[dr] \ar@{-->}[rr]^{\sigma_R} & & \bullet \ar[ur] \ar[dr] & \hdots \\
  \bullet \ar[ur] \ar[r]  \ar[dr] \ar@{-->} @/^.7pc/ [rr]^{\sigma_R} & \bullet \ar[r]  \ar@{-->}@/_.5pc/[drr]_(.8){\sigma_R} & \bullet \ar[ur] \ar[r] \ar[dr] \ar@{-->}@/^.7pc/[rr]^{\sigma_R} & \bullet \ar[r]  \ar@{-->}@/_.5pc/[drr]_(.8){\sigma_R} & \bullet \ar[ur] \ar[r] \ar[dr] \ar@{-->}@/^.7pc/[rr]^{\sigma_R} & \bullet \ar[r] & \bullet \\
   & \bullet \ar[ur]  \ar@{-->}@/_.5pc/[urr]_(.2){\sigma_R} &  & \bullet \ar[ur]  \ar@{-->}@/_.5pc/[urr]_(.2){\sigma_R} &  & \bullet \ar[ur]}$$

Thus we have $\Omega^4\sigma_R^2(k) \cong k$ for the unique simple $A$-module $k$.  By Proposition~\ref{prop:LinckelmannTheorem}, we can conclude that $\sigma_R^2 \cong \Omega^{-4} \circ F$ on $\stmod A$, where $F$ is a Morita auto-equivalence of $\rmod A$.  Moreover $F$ induces the identity automorphism on this component of the stable AR-quiver.


\subsection{Extraspecial $p$-groups}  The construction of the $\mathbb{P}^n$ stable twists in the last section was partially motivated by Alperin's construction of endo-trivial modules using relative syzygies \cite{Alp}, as described in Theorem 3.1 of \cite{CarThe2}.  In particular, for an odd prime $p$, we focus on an extraspecial $p$-group $G$ of order $p^3$ and exponent $p$, and demonstrate that $\mathbb{P}^{p-1}$ stable twists are defined for $kG$ (with $k$ an algebraically closed field of characteristic $p$), and that the endo-trivial $kG$-modules are recovered as the images of the trivial module $k$ (and its syzygies) under these twists.

Let $G = \gen{g, h, z\ |\ g^p=h^p=z^p=1, z=[g,h], gz = zg, hz=zh}$ be the extraspecial $p$-group of order $p^3$ and exponent $p$ (we follow the notation of \cite{Car2}, \S 6).  Let $H_i = \gen{g^ih}$ and $E_i = \gen{z,g^i h} = C_G(H_i)$ for $1 \leq i \leq p$, and set $H_{p+1} = \gen{g}$ and $E_{p+1} = \gen{z,g} = C_G(H_{p+1})$.  We also let $A = kG$ and $R_i = kH_i$ for each $1 \leq i \leq p+1$.  Notice that $x_i := 1-g^ih$ (resp. $x_{p+1} = 1-g$) generates $kH_i$ as an algebra and we make the identification $kH_i = k[x_i]/(x_i^p)$.  Furthermore $y = 1-z \in kG$ commutes with each $x_i$.  We set $T_i = A/x_i A = k_{H_i}\uparrow^G$.  In order to show that $\rho =\rho_{R_i,y}$  actually defines a $\mathbb{P}^{p-1}$ stable twist, it suffices to verify $\stEndo_{kG}(T_i) \cong k[\psi]/(\psi^p)$ with $\psi$ corresponding to multiplication by $y$.    

To simplify notation, we fix $i=p$ so that we can work with $H=H_p$ and $R = R_p$.  (In fact, as $G$ has automorphisms sending $g \mapsto g, h \mapsto g^ih, z \mapsto z$ or sending $g \mapsto h, h \mapsto g, z \mapsto z^{-1}$, the other cases can be deduced from our calculations here.  Moreover, it will follow that all the $\rho_{R_i,y}$ differ from one another by Morita equivalences.)  With respect to the subgroup $H$, $G$ has the double coset decomposition $$G = \bigcup_{j=0}^{p-1} Hz^jH \cup \bigcup_{l=1}^{p-1} H g^l H,$$ with $Hz^jH = Hz^j$ as $z$ is central.  It follows that as a $(kH,kH)$-bimodule $$kG = \bigoplus_{j=0}^{p-1} k(Hz^jH) \oplus \bigoplus_{l=1}^{p-1} k(H g^l H) \cong \bigoplus_{j=0}^{p-1} (kH \otimes z^j) \otimes_{kH} kH \ \oplus \ \bigoplus_{l=1}^{p-1} (kH \otimes g^l) \otimes_k kH,$$
as $H \cap H^{z^j} = H$ for each $j$ and $H \cap H^{g^l} = 1$ for each $l$.  Therefore, by Mackey's theorem we have
$$k_H \uparrow^G \downarrow_{H} \cong \bigoplus_{j=0}^{p-1} (k_H \otimes z^j) \oplus \bigoplus_{l=1}^{p-1} (k \otimes g^l)\uparrow^H \cong k_H^p \oplus kH^{p-1},$$
where $k_H \otimes z^j \cong k_H$ for each $j$ and $(k \otimes g^l)\uparrow^{H} \cong kH$ for each $l$. In particular, for $T = k_H \uparrow^G$, we see that $\stHom_{kG}(T,T) \cong \stHom_{kH}(k_H,T\downarrow_H)$ is $p$-dimensional.  Moreover, the automorphism of $T$ induced by multiplication by $z$ restricts to an automorphism of the non-projective part of $T\downarrow_H$, cyclicly permuting the $p$ components of the form $k \otimes z^j$.  Let $\psi$ be the endomorphism of $T$ induced by multiplication by $y=1-z$.  Clearly $\psi^p=0$ and the restriction of $\psi^j$ to an endomorphism of the non-projective part of $T \downarrow_H$ is non-zero for each $j<p$.  In fact, this restriction is non-zero in the stable category as well since it is an endomorphism of a (non-projective) semisimple module.  It follows that the non-zero powers of $\psi$ are linearly independent in $\stEndo_{kH}(T\downarrow_H)$ and hence in $\stEndo_{kG}(T)$ as well.  Hence we conclude that $\stEndo_{kG}(T) = k[\psi]/(\psi^p)$ as required.

Finally, Proposition 6.5 shows that the image of the trivial module $k_G$ under $\rho_{R_i,y}$ coincides with the endo-trivial module $N_i$ described in Theorem 3.1(c) of \cite{CarThe2}.  According to Theorem 6.1 of \cite{CarThe2}, these endo-trivial modules for $2 \leq i \leq p+1$, together with $\Omega(k)$, generate the group of endo-trivial $G$-modules, which is isomorphic to $\mathbb{Z}^{p+1}$.

\subsection{Discussion and open questions}

Naturally, one would like to understand how general our constructions are, and find more examples of local symmetric algebras where they yield non-trivial auto-equivalences of the stable module category.  While finding examples of algebras $A$ and $R$ satisfying the hypotheses of Section 4 is not difficult, the condition that the stable endomorphism ring $\stEndo(k\otimes_R A)$ is isomorphic to $k[\psi]/(\psi^{n+1})$ seems harder to satisfy.  We also do not know of any examples where $m \neq n+1$ in our notation (recall $R \cong k[t]/(t^m)$ and $\stEndo_A(k\otimes_R A) \cong k[\psi]/(\psi^{n+1})$), although from our proofs we see no immediate reason why these numbers would need to coincide.  It would also be interesting to see other examples of non-trivial auto-equivalences of stable module categories of local symmetric algebras which do not belong to the group of auto-equivalences generated by those constructed here and the syzygy functor.  
 
 One important class of examples is furnished by the group algebras of the generalized quaternion $2$-groups in characteristic two.  Such group algebras admit an endo-trivial module that is not a syzygy of $k$, and we note that neither spherical stable twists nor $\mathbb{P}^n$ stable twists appear to recover the auto-equivalences induced by this endo-trivial module.  Based on one construction of these endo-trivial modules found in Theorem 3.1(c) of \cite{CarThe2}, it seems plausible that the bimodule inducing these auto-equivalences may be obtained by a {\it quadruple} cone construction, which could then be generalized to arbitrary local symmetric algebras of quaternion type.  Moreover, we expect that the occurence of a quadruple cone here should be tied to the fact that the local symmetric algebras of quaternion type are periodic algebras of period $4$ \cite{ErdSko}, meaning that they have periodic projective resolutions as modules over their enveloping algebras.
 
 Such a connection to periodicity is further supported by Grant's recent work on periodic twists \cite{Gra}.  For starters, in the derived category of coherent sheaves on a projective variety $X$, spherical twists and $\mathbb{P}^n$ twists are associated to objects with graded endomorphism rings isomorphic to $k[x]/(x^2)$ and $k[x]/(x^{n+1})$ respectively; both of which are periodic algebras.    Grant then proceeds to define auto-equivalences of the derived category of a symmetric $k$-algebra $\Lambda$ associated to any projective $\Lambda$-module with a periodic endomorphism algebra, and shows that spherical twists and $\mathbb{P}^n$ twists (in this context) can be recovered as special cases.  Moreover, it becomes clear that the single and double cone constructions involved in the definition of spherical and $\mathbb{P}^n$ twists are tied to the $1$ and $2$-periodicity, respectively, of the algebras $k[x]/(x^2)$ and $k[x]/(x^{n+1})$.  As we enconter in our constructions the same two endomorphism rings of objects in the stable category $\stmod A$, it motivates us to ask if one can define additional auto-equivalences of $\stmod A$ for a symmetric algebra $A$ tied to objects in $\stmod A$ with periodic endomorphism algebras.  More generally still, one wonders whether there is a single construction of an exact auto-equivalence of a triangulated $k$-category associated to an object with a periodic endomophism ring, which unifies all of these constructions.

\end{document}